\magnification=\magstep1%metka
\input amstex
\documentstyle{amsppt}
\TagsOnRight
\hsize=5in                                                  
\vsize=7.8in
\strut
\vskip5truemm
\font\small=cmr8
\font\itsmall=cmti8
%\font\ref=cmr9
%\font\refit=cmti9
%\font\refbf=cmbx9

\def\smallarea#1{\par\begingroup\baselineskip=10pt#1\endgroup
\par}
\def\abstract#1{\begingroup\leftskip=5mm \rightskip 5mm
\baselineskip=10pt\par\small#1\par\endgroup}

\def\refer#1#2{\par\begingroup\baselineskip=11pt\noindent
\leftskip=8.25mm\rightskip=0mm

\strut\llap{#1\kern 1em}{#2\hfill}\par\endgroup}
              
\nopagenumbers
%\subheading{Revised version}
%\vskip 2cm
\centerline{\bf HOMOLOGICAL ALGEBRA OF NOVIKOV-SHUBIN INVARIANTS} 
\centerline{\bf AND MORSE INEQUALITIES}

%\vskip 25.3pt
\vskip 2cm

\centerline{MICHAEL S. FARBER
\footnotemark"$^1$"}
%{${}^*$}{\baselineskip=12pt\smalltype 
\footnotetext"$^{1}$"{The research was supported
by a grant from US - Israel binational Science Foundation}
\smallarea{%
\centerline {\itsmall School of Mathematical Sciences,}
\centerline {\itsmall Tel-Aviv University,}
\centerline {\itsmall Tel-Aviv 69978, Israel}
\centerline {\itsmall e-mail farber$\@$math.tau.ac.il}
\centerline {\itsmall fax 972-3-6407543}
}
\medskip

\define\C{{\Bbb C}}
\define\R{{\Bbb R}}  
       
\define\Z{{\Bbb Z}}

\define\T{{\Cal T}}  
\define\X{{\Cal X}}   
\define\Y{{\Cal Y}}  

\define\N{{\Bbb N}}
 
\define\Hom{\operatorname{Hom}}

\define\Tor{\operatorname{Tor}} 
\define\TOR{{\Cal Tor}}
\define\EXT{{\Cal Ext}}
\define\Tot{\operatorname{Tot}}
   
\define\Ext{\operatorname{Ext}}

\define\Ker{\operatorname{Ker}} 
\define\End{\operatorname{End}}

\define\A{{\Cal A}}  
\redefine\B{{\Cal B}}  
\define\e{{\frak e}} 
\redefine\c{{\frak c}}

\define\im{\operatorname{im}} 
\redefine\Im{\operatorname{Im}}   
%\define\Ker{\operatorname{Ker}}   
\define\Coker{\operatorname{Coker}}   
\define\Coim{\operatorname{Coim}}   
\define\id{\operatorname{id}}    
\define\ns{\operatorname{ns}}     
\define\ob{\operatorname{ob}}   
\define\cl{\operatorname{cl}} 
%\define\ln{\operatorname{ln}}
\redefine\H{\Cal H}

\define\E{\Cal E}
  
\def\<{\langle}
\def\>{\rangle}

\define\pd#1#2{\dfrac{\partial#1}{\partial#2}}
\documentstyle{amsppt}   

\nopagenumbers

\vskip 2cm
\abstract{\subheading{Abstract} It is shown in this paper that the 
topological phenomenon {\it "zero in the continuous spectrum"}, 
discovered by S.P.Novikov and 
M.A.Shubin, can be explained in terms of a homology theory on the category
of finite polyhedra with values 
in certain abelian category. This approach implies homotopy invariance of the
Novikov-Shubin invariants. Its main advantage is that it allows to use the 
standard homological techniques, such as spectral sequences, derived functors,
universal coefficients etc., while studying the Novikov-Shubin invariants. 
It also leads to some new quantitative invariants, measuring the 
Novikov-Shubin phenomenon in a different way, which are used in the present
paper in order to strengthen the Morse type inequalities 
of S.P. Novikov and M.A. Shubin \cite{NS1}.}
\vskip 0.5cm

\heading{\bf \S 0. Introduction}\endheading

This paper suggests a conceptually new approach, which unites the
$L^2$ cohomology theory and the Novikov-Shubin invariants. 
It is shown here that these theories are two 
different parts of {\it a unique cohomology theory with values in an 
abelian category}. This abelian category, denoted $\E(\A)$, 
contains the familiar additive 
category of Hilbertian modules over a von Neumann algebra $\A$ as a full
subcategory of projectives.
An important abelian subcategory
of $\E(\A)$ is formed by {\it torsion virtual Hilbertian modules}. It turns 
out that any object of $\E(\A)$ has 
canonically defined {\it torsion and projective parts} and coincides 
with their direct sum. The von Neumann dimension 
is an invariant of the projective part; similarly,
the Novikov-Shubin number is an invariant of the torsion part.

There are natural homology and cohomology theories 
with values in the
abelian category $\E(\A)$. I denote these theories by $\H_i(X,M)$ and 
$\H^i(X,M)$ correspondingly and call {\it extended $L^2$ homology and 
cohomology}. Here $X$
is a CW complex having finitely many cells in every dimension, 
and $M$ is a Hilbertian $(\A-\pi)$- or $(\pi-\A)$-bimodule
(cf. \S 6 below), and $\pi=\pi_1(X)$ is the fundamental group of $X$.  
These theories are homotopy invariant.
The projective part of $\H^i(X,M)$ coincides
with the usual {\it reduced} $L^2$ cohomology. The Novikov-Shubin invariants
of the torsion part of $\H^i(X,M)$ coincide with the Novikov-Shubin 
invariants of the complex $X$. This gives a conceptually transparent proof 
of homotopy invariance of both von Neumann Betti numbers (proven by J.Dodziuk 
\cite{D}) and of
the Novikov-Shubin invariants (proven by M.Gromov and M.Shubin \cite{GS}).

Simple examples of torsion virtual Hilbertian modules show that they can be 
not isomorphic but have equal Novikov-Shubin invariants.

I introduce here a new numerical invariant of torsion objects, which is 
independent on the Novikov-Shubin number; it is called {\it the minimal
number of generators}.

As the main application, new Morse type inequalities for the 
numbers of critical points of a function on compact manifold are found. 
These 
inequalities give quantitative information even in the cases when all
von Neumann Betti numbers vanish.
Note, that Morse inequalities using the von Neumann dimensions of the reduced
$L^2$ cohomology were obtained by S.P.Novikov and M.A.Shubin in \cite{NS}.
A more complete exposition will appear in \cite{Sh}, where a general 
philosophy of {\it model operator} is developed.

The most important advantage of the suggested in this paper new 
approach to the Novikov-Shubin invariants consists, perhaps, in the fact that
this approach allows to apply the standard techniques of homological algebra 
(such as spectral sequences, derived functors, etc.) to studying the
Novikov-Shubin invariants. Some results in this direction are included in this
paper. For example, a universal coefficients theorem, and a Poincar\'e 
duality theorem, are proven here. We also study the group homology and 
cohomology with values in the extended abelian category $\E(\A)$. It is
shown that
in the most general case the extended $L^2$ homology and cohomology can be
expressed through these group cohomology by means of some Cartan-Eilenberg
type spectral sequences; the initial $E^2$-term of these spectral sequences 
depends only on the homology of the universal covering, considered as
modules over the group ring. This allows to 
conclude, for example, that the Novikov-Shubin invariants of a space with
free fundamental group depend only on the homology modules of the universal
covering.

Before the present work has started, M.A.Shubin had shown to me a simple 
argument (due to M.Gromov) which proves that nontriviality of the 
Novikov-Shubin number implies existence of at least one critical point. 
M.A.Shubin also suggested to me the problem of finding quantitative 
estimates. This problem was the main motivation for the present work.
I am also thankful to W.L\"uck who sent to me preprint of \cite{LL} 
before publication; it was very influential. Finally,
I would like to mention interesting and very helpful discussions of 
different parts of this work with M.Gromov,
J.Levine, V.Mathai and S.Weinberger. 

This paper distributed as a preprint of the Tel Aviv University in July 1995; 
a short announcement was published in \cite{Fa}. 

Wolfgang L\"uck, in his recent preprint "Hilbert modules and modules over
finite von Neumann algebras and applications to $L^2$-invariants" (Mainz, 
December 1995) suggested a different (more algebraic) approach to studying
$L^2$-invariants. He also discussed some interesting examples and applications, 
and showed that his approach is in fact equivalent to the method of this paper.

\heading{\bf \S 1. The category of Hilbertian modules over a von Neumann 
algebra}\endheading

This section briefly describes a modification of the standard category of 
Hilbert modules over a finite von Neumann algebra. A more
detailed exposition of this material can be found in \cite{CFM}.

\subheading{1.1} Let ${\A}$ be a finite  von Neumann algebra 
with a fixed finite,
normal, and faithful trace $\tau:{\Cal A}\rightarrow \C$.  The involution
in ${\A}$ will be denoted $*\,$; by $\ell^{2}({\A})$ we denote the
completion of ${\A}$ with respect to the scalar product $\langle a,b
\rangle=\tau(b^{*}a)$, for $a,b\in{\Cal A}$.

Recall that a {\it Hilbert module} over ${\A}$ 
is a Hilbert space $M$
together with a continuous left ${\Cal A}$-module structure such that there
exists an isometric ${\Cal A}$-linear embedding of $M$ into 
$\ell^{2}({\Cal A})\otimes H$, for some Hilbert space $H$.  
Note that this embedding is not
part of the structure.  A Hilbert module $M$ is {\it finitely generated} if it
admits an embedding $M\rightarrow\ell^{2}({\Cal A})\otimes H$ as above with
finite dimensional $H$.

Any Hilbert module, being a Hilbert space, has a particular scalar product.  
In this paper we
wish to consider a weaker notion obtained from Hilbert module by forgetting
the scalar product but preserving its topology and the ${\Cal A}$-action.

\subheading{1.2. Definition} A {\it  Hilbertian module} 
is a topological vector space $M$
with continuous left ${\Cal A}$-action such that there exists a scalar product
$\langle\;,\;\rangle$ on $M$ which generates the topology of $M$ and such
that $M$ together 
with $\langle\;,\;\rangle$ and with the ${\Cal A}$-action is a Hilbert
module. In particular, the involution on $\A$ is compatible with the 
involution on the space of bounded linear operators on $M$ determined by the
scalar product $\<\ ,\ \>$.

If $M$ is a Hilbertian module, 
then any scalar product $\langle\;,\;\rangle$ on
$M$ with the above properties will be called {\it admissible}.  

\subheading{1.3} It is easy to see that {\it given a Hibertian module, 
all different choices of
admissible scalar products on it produce isomorphic Hilbert modules}, 
cf. \cite{CFM}. The situation here is
similar to the case of finite dimensional vector spaces: any
choice of a scalar product on a vector space
produces an isomorphic Euclidean vector space.

Using the above mentioned fact, we may define 
{\it finitely generated Hilbertian modules}
as those for which the corresponding Hilbert modules (obtained by a 
choice of an admissible scalar product) are finitely generated,
cf. 1.1 above.

Similarly, we may correctly define {\it the von Neumann dimension
$\dim_\tau(M)$ of a Hilbertian module $M$} as the von Neumann dimension
of the Hilbert module obtained by a choice of an admissible scalar product
on $M$.

\subheading{1.4} Let us denote by $\H(\A)$ the category, whose objects
are finitely generated Hilbertian modules over $\A$ and whose morphisms are
continuous linear maps commuting with the action of the algebra $\A$.
Obviously, $\H(\A)$ is an additive category.

Note that the category $\H(\A)$ depends on the choice of the trace
$\tau$ in an essential way, although the trace $\tau$ does not appear in 
the notation $\H(\A)$.

Given a morphism $f:M\to N$ in the category $\H(\A)$, the set-theoretic
kernel $\ker(f)\subset M$ of $f$ has naturally a structure of finitely 
generated Hilbertian module over $\A$; it coincides 
with the kernel of $f$ in sense of the category theory.
On the contrary, the set-theoretic image 
$\im(f)\subset N$ is not in general closed and so it is not a Hilbertian
module; the categorical image of the morphism $f$ coincides with 
the closure of the set-theoretic image $\cl(\im(f))\subset N$.

Let's emphasize that the symbols $\ \im$ and $\ \ker$ will always denote the
set-theoretic notions. The symbols $\Im, \ \Ker,\ \Coim, \ \Coker$
(starting with the capital letters)
will denote the corresponding notions of the category theory, 
cf., for example, \cite{Gr}, \cite{F}.

\subheading{1.5} Note the following well-known important property of 
the category $\H(\A)$. If 
$$0\to M^\prime@>{f}>>M@>g>>M^{\prime\prime}\to 0$$
is a sequence of objects and morphisms of the category $\H(\A)$, which is 
exact, considered as a sequence of  
abelian groups (i.e. if $f$ is injective, $g$ is
surjective, and $\im(f)=\ker(g)$), then this sequence
{\it splits in the category $\H(\A)$}.

\subheading{1.6. Duality} There is a notion of duality in $\H(\A)$.
Given a Hilbertian module $M$, consider the set $M^\ast$ of all anti-linear
continuous functionals on $M$. Introduce the following action of $\A$ on 
$M^\ast$: if $\phi\in M^\ast$ and $\lambda\in \A$ then 
$(\lambda\cdot \phi)(m)=\phi(\lambda^\ast\cdot m)$ for all $m\in M$. 
Here $\lambda^\ast$ is defined by the involution $\ast$ of $\A$. 
In particular, this defines an action of $\C\subset \A$ on $M^\ast$. 

Note that, the dual of
$\ell^2(\A)$ is canonically isomorphic to $\ell^2(\A)$. If $M$ is an arbitrary
Hilbertian module then $M^\ast$ is also a Hilbertian module, which is 
isomorphic to $M$, but not canonically. 
Any choice of an admissible
scalar product on $M$ gives an isomorphism between $M$ and $M^\ast$. 
%Conversely, any isomorphism between Hilbertian modules $M$ and $M^\ast$
%determines an admissible scalar product on $M$.
Note that $M$ is finitely generated if and only if $M^\ast$ is. 

Duality is a contravariant functor: if $f:M\to N$ is a morphism of Hilbertian
modules then the dual map $f^\ast :N^\ast\to M^\ast$ is also a morphism
of Hilbertian modules. Duality is also involutive:  
$M^{\ast\ast}$ is canonically isomorphic to $M$. The canonical isomorphism 
$M\to M^{\ast\ast}$ is given by $m\mapsto(f\mapsto\overline{f(m)})$, where
$m\in M$, and $f\in M^\ast$. Here the bar denotes the complex conjugation.

\heading{\bf \S 2. The Extended abelian category}\endheading

In this section an abelian category, containing the category of Hilbertian
modules $\H(\A)$, is described. 

$\A$ denotes in this section a finite von Neumann algebra
supplied with a fixed finite, normal, and faithful trace $\tau$.

\subheading{2.1} The category $\H(\A)$ of Hilbertian modules 
over the von Neumann algebra $\A$ {\it is not an abelian category}.
In fact, the condition AB1 of \cite{Gr} is satisfied (any morphism has a 
kernel and a cokernel) but the condition AB2 in general does not hold. 
Recall that the 
condition AB2 requires that the canonical map $\Coim(f)\to \Im(f)$ be
an isomorphism. Given a morphism $f:M\to N$ in the category $\H(\A)$,
its coimage $\Coim(f)$ can be represented by the Hilbertian module
$M/\ker(f)$, while (as it was mentioned above in 1.4) $\Im(f)=\cl(\im(f))$.
Thus, the canonical map $\Coim(f)\to \Im(f)$ is a continuous linear map,
which is injective and the set of its values is dense. Such map clearly 
may be not invertible in $\H(\A)$; one easily constructs the corresponding
examples.

A few possible approaches to the general problem of {\it representing} a 
given additive category in abelian
categories were discussed by P.Freyd in \cite{F}. One of the constructions,
described in \cite{F}, gives the following abelian category $\E(\A)$,
which we will call {\it the extended category of Hilbertian modules}.
Note that the category $\E(\A)$ depends on the choice of the trace $\tau$
although $\tau$ does not appear explicitly in the notation.

\subheading{2.2. Definition} An {\it object} of the category $\E(\A,\tau)$ 
is defined as a
morphism $(\alpha:A^\prime\to A)$ in the category $\H(\A)$. 
Recall that here $A^\prime$ and
$A$ are finitely generated Hilbertian modules over the von Neumann algebra
$\A$ and $\alpha$ is a continuous
linear map commuting with the action of the algebra $\A$.

Given a pair of objects $\X=(\alpha: A^\prime \to A)$ and 
$\Y=(\beta:B^\prime\to B)$ of $\E(\A)$,
a {\it morphism} $\X\to\Y$ in the category $\E(\A)$ is an  
equivalence class of morphisms $f:A\to B$ of category $\H(\A)$
such that $f\circ\alpha=\beta\circ g$
for some morphism $g:A^\prime \to B^\prime$ in $\H(\A)$. 
Two morphisms $f:A\to B$ and $f^\prime:A\to B$ of $\H(\A)$ represent 
{\it identical morphisms of $\E(\A)$}
$\X\to\Y$ iff 
$f-f^\prime = \beta\circ F$ for some morphism $F:A\to B^\prime$ of category
$\H(\A)$. This defines an equivalence relation.
The morphism $\X\to\Y$, represented by $f:A\to B$
well be denoted 
$$[f]:(\alpha:A^\prime\to A)\ \to\ (\beta:B^\prime\to B)\quad\text{or}\quad
[f]:\X\to\Y.$$
Sometimes we will say that $[f]$ is represented by a diagram
$$
\CD
(A^\prime@>\alpha>>A)\\
 & &                @VVfV\\
(B^\prime@>\beta>>B)
\endCD\tag1
$$
this emphisizes existence of the morphism $g:A^\prime\to B^\prime$ making 
the diagram commutative.

The {\it composition} of morphisms is defined as the composition of the 
corresponding morphisms $f$ in the category $\H(\A)$; this unambiguously 
defines a composition law for morphisms. 
The {\it identity} morphism is $[\id_A]$.
Clearly, $\E(\A)$ is an additive category.

Objects of $\E(\A)$ will be called {\it virtual Hilbertian modules}.

It follows from the work of P.Freyd \cite{F}, that
{\it the extended category $\E(\A)$ is an abelian category}. 

Having in mind our further purposes in this paper,
we are going to sketch an independent proof. 

First, we have to be able 
to compute kernels and cokernels of morphisms in $\E(\A)$.

\proclaim{2.3. Proposition} Suppose that
$[f]:(\alpha:A^\prime\to A)\ \to\ (\beta:B^\prime\to B)$          
is a morphism in the extended category of Hilbertian modules $\E(\A)$.
Then its kernel is represented by
$$[k]:\ (\gamma:P^\prime \to P)\ \to\ (\alpha:A^\prime\to A),$$
where
$$
\botaligned
{\CD
P@>k>>A\\
@VVV  @VVfV\\ 
B^\prime@>\beta>>B&
\endCD}
\qquad\text{{and}}\qquad
{\CD
P^\prime@>{k^\prime}>>A^\prime\\
@VVV  @VV{f\circ \alpha}V\\ 
B^\prime@>\beta>>B&
\endCD}
\endbotaligned
$$
are the pullbacks of the diagrams

$$
\botaligned
{\CD
  @. A\\
@.   @VVfV\\ 
B^\prime@>\beta>>B
\endCD}
{\qquad\text{and}\qquad}
{\CD
  @. A^\prime\\
@.   @VV{f\circ\alpha}V\\ 
B^\prime@>\beta>>B
\endCD}
\endbotaligned\tag2
$$
correspondingly, and $\gamma:P^\prime\to P$ is the canonical map, induced 
by the obvious map of the right diagram (2) into the left one.

The cokernel of the above morphism $[f]$
%$$[f]:(\alpha:A^\prime\to A)\ \to\ (\beta:B^\prime\to B),$$  
is represented by
$$[\id_B]:\ (\beta:B^\prime\to B)\ \to\ ((\beta,-f):B^\prime\oplus A\to B).$$
\endproclaim

\demo{Proof} It follows by checking straightforwardly the definitions.
We will leave this proof to the reader. $\square$
\enddemo

Thus, we obtain that {\it any morphism of $\E(\A)$ has a kernel and a 
cokernel} and so condition AB1 of \cite{Gr}, \S 1.4 is satisfied.

Using Proposition 2.3, we may compute explicitly 
the coimage $\Coim([f])$
(which is defined as the cokernel of the kernel) and the image 
$\Im([f])$ (which is defined as the kernel of the cokernel); 
the result in both cases coincides with
$$[f]:(k:P\to A)\to (\beta:B^\prime\to B)$$
where $P$ and $k$ are the same as in Proposition 2.3.

This shows that the condition AB2 from the definition of abelian categories,
cf. \cite{Gr}, chapter 1, is satisfied. {\it Hence $\E(\A)$ is an abelian
category}.

\subheading{2.4. Excision} Note that two very different maps
$(\alpha:A^\prime\to A)$ and $(\beta:B^\prime\to B)$ may represent isomorphic
objects of category $\E(\A)$. In order to clarify this question,
consider the following situation.
 
Suppose that $(\alpha:A^\prime\to A)$ is an object of $\E(\A)$ and
$P\subset A^\prime$ is a closed $\A$-submodule such that its image
$\alpha(P)\subset A$ is also closed. Let $B^\prime$ and $B$ denote the
factor-modules $B^\prime\ =\ A^\prime/P$ and $B\ =\ A/\alpha(P)$.
The map $\alpha$ induces the obvious map $\beta:B^\prime\to B$. We
claim now that {\it the obtained object $(\beta:B^\prime\to B)$ of $\E(\A)$ 
is isomorphic to the initial $(\alpha:A^\prime\to A)$}. 

The passage from 
$(\alpha:A^\prime\to A)$ to $(\beta:B^\prime\to B)$ described above,
will be called {\it excision with respect to $P$}. We will also say that
$(\alpha:A^\prime\to A)$ is an {\it enlargement of} $(\beta:B^\prime\to B)$.

Observe first that our claim is immediate in the case $\alpha(P)=0$. 

It follows that we may always make an excision with respect to the kernel
of the map $\alpha$ and so 
{\it any object of $\E(\A)$
can be represented by $(\alpha:A^\prime\to A)$  
with an injective morphism $\alpha$}.

Suppose now that $(\beta:B^\prime\to B)$ is obtained from
$(\alpha:A^\prime\to A)$ with injective $\alpha$
by excision with respect to $P\subset A$. Then the following sequence
in $\H(\A)$
$$0\to A^\prime@>{
\bmatrix f^\prime\\ \alpha\endbmatrix}>>
B^\prime\oplus A@>{(-\beta,f)}>>B\to 0$$
is exact. Here 
$f^\prime:A^\prime\to A^\prime/P\ =\ B^\prime$ and 
$f:A\to A/P\ =\ B$ are the canonical projections. Hence this sequence 
splits (by 1.5). Consider a splitting
of the above sequence
$$0@<<< A^\prime@<{(\delta^\prime,\sigma^\prime)}<<
B^\prime\oplus A@<{
\bmatrix \sigma\\ \delta\endbmatrix}<<
B@<<< 0$$
Then the relations
$$f\circ\alpha\ =\ \beta\circ f^\prime\qquad\text{and}\qquad
\delta\circ\beta\ =\ \alpha\circ\delta^\prime$$
show that $[f]:(\alpha:A^\prime\to A)\to(\beta:B^\prime\to B)$ and
$[\delta]: (\beta:B^\prime\to B)\to (\alpha:A^\prime\to A)$ are morphisms of
$\E(\A)$, and the relations
$$-\beta\circ\sigma+f\circ\delta\ =\ 1_B\qquad\text{and}\qquad
\alpha\circ\sigma^\prime + \delta\circ f = 1_A$$
show that the above morphisms $[f]$ and $[\delta]$ are mutually inverse.
This completes the proof.

Note that {\it a morphism $(\alpha:A^\prime\to A)$ represents a null object
of $\E(\A)$} (which can be characterized as an object with the property that 
its identity morphism coincides with the zero morphism)
{\it if and only if $\alpha$ is surjective}.

\subheading{2.5. Monomorphisms and epimorphisms} Let 
$[f]:(\alpha:A^\prime\to A)\ \to\ (\beta:B^\prime\to B)$ 
be a morphism of $\E(\A)$. It is easy to see (using Proposition 2.3)
that $[f]$ {\it is a monomorphism of the category $\E(\A)$ if and only if }
$$\alpha(A^\prime)\supset f^{-1}(\beta(B^\prime))\tag3$$
and {$[f]$ {\it is an epimorphism of the category $\E(\A)$ if and only if} 
$$B\ = \ \beta(B^\prime)\ +\ f(A)\tag4$$

In particular, any morphism $f:A\to B$ of $\H(\A)$ determines a morphism
$[f]:(0\to A)\to (0\to B)$ in $\E(\A)$ and $f$ is injective (respectively,
surjective) as a morphism vector spaces, if and only if the morphism
$[f]$ is a monomorphism (respectively, epimorphism).

The following remark will be useful later.
\proclaim {2.6. Lemma} Given a monomorphism in $\E(\A)$
$$[f]:(\alpha:A^\prime\to A)\ \to\ (\beta:B^\prime\to B),$$ 
one can perform an excision on $(\alpha:A^\prime\to A)$
such that the same monomorphism will be represented by a diagram
$$
\CD
(C^\prime@>\gamma>>C)\\
 & &                @VVgV\\
(B^\prime@>\beta>>B)
\endCD\tag5
$$
with injective $g$.
\endproclaim
\demo{Proof} Let $P^\prime$ denote $\ker f\subset A$ and let 
$P=\alpha^{-1}(P^\prime)$. Then $\alpha(P)=P^\prime$ (by virtue of (3))
is closed 
and so we can perform an excision with respect to $P$.
This completes the proof.
$\square$\enddemo

We can now strengthen Lemma 2.6:
\proclaim{2.7. Corollary} (1) For any monomorphism $\X\to\Y$ in $\E(\A)$, 
one can perform a sequence of excisions and enlargements on $\X$ and on $\Y$
such that the given monomorphism will be represented
by a diagram
$$
\CD
(A^\prime@>\alpha>>A)\\
 @VhVV   @VVfV\\
(B^\prime@>\beta>>B)
\endCD\tag6
$$
with $\alpha$,$\beta$,$f$ monomorphisms and $h$ isomorphism. (Note, that
$h$ is determined uniquely in this situation by $\alpha$, $\beta$, $f$.)
Conversely,
any morphism of $\E(\A)$ represented by such diagram is a monomorphism.

(2) For any epimorphism $\X \to \Y$ in $\E(\A)$ one can perform a sequence of 
excisions and enlargements on $\X$ and on $\Y$ such that the given epimorphism
will be represented by a diagram (6) as above
with $\alpha$,$\beta$, and $h$ injective and $f$ bijective. (Here again
$h$ is determined uniquely by $\alpha$, $\beta$ and $f$.) Conversely, any
such diagram represents an epimorphism in $\E(\A)$. 
\endproclaim
\demo{Proof} (1) First we represent a given monomorphism by a diagram (6)
with $\beta$ injective and $h$ surjective; this can easily be arranged. 
Then we perform excision with respect
to $\ker \alpha$ and then with respect to $\ker(f)$ (as explained in 2.6).
This will give us a diagram representing the given morphism having the 
desired properties. The converse statement follows from 2.5. 

(2)  First, represent a given epimorphism by a diagram of form (6) 
with $\beta$ injective and $f$ surjective. This can be done starting from 
an arbitrary representation with injective $\beta$ by considering the diagram
$$
\CD
A^\prime\oplus B^\prime @>{\bmatrix \alpha & 0\\ 0 & 1\endbmatrix}>>
A\oplus B^\prime\\
@V{[h,1]}VV      @VV{[f,\beta]}V\\
B^\prime @>>{\beta}>B.
\endCD
$$
Assuming that this has been arranged, we perform a similar enlargement such
that the image of the new map $\alpha$ would contain the kernel of $f$.
Now consider $\ker(f)=P^\prime$ and $\alpha^{-1}(P^\prime)=P$. Performing
excision with respect to $P$ (using the arguments similar to those 
of Lemma 2.6) we obtain a diagram representing the given epimorphism which
has $\beta$ injective, $f$ bijective and $\alpha$ injective. Then $h$ is 
injective, and the result follows. The converse statement follows obviously
from 2.5. $\square$\enddemo

\subheading{2.8. Embedding of $\H(\A)$ into $\E(\A)$} Given a Hilbertian
module $A$, consider the zero morphism $(0\to A)$ as an object of $\E(\A)$.
Any morphism $f:A\to B$ in $\H(\A)$ determines the morphism 
$[f]: (0\to A)\to (0\to B)$
in $\E(\A)$;
conversely, any morphism in $\E(\A)$ 
$$(0\to A)\to (0\to B)$$
determines uniquely a morphism of $\H(\A)$ between $A$ and $B$. This shows 
that we have a {\it functor 
$$\frak F:\H(\A)\ \to \E(\A),$$
and this functor is a full embedding}.

Let us show that {\it an object $\X\in\ob(\E(\A))$ is isomorphic to a 
Hilbertian module in $\E(\A)$ if and only if $\X$ is projective.}

Note that it is not true that any projective $\X\in\ob(\E(\A))$ is a
Hilbertian module, i.e. comes from $\H(\A)$.

First we will show that any Hilbertian module $(0\to X)$ is projective
in $\E(\A)$. Suppose that we have a diagram 
$$
\CD
(A^\prime@>\alpha>>A)\\
 & &                @VVfV\\
(B^\prime@>\beta>>B)\\
@.                @AAgA\\
(0@>>>X)
\endCD\tag7
$$
with morphism $[f]$ being epimorphism (as a morphism of $\E(\A)$); this 
means that $B=\beta(B^\prime)+f(A)$ by 2.5. Then we have an exact sequence
in $\H(\A)$
$$B^\prime\oplus A@>{(-\beta, f)}>>B\to 0,$$
and so the map $g:X\to B$ can be lifted into $B^\prime\oplus A$ (by 1.5), 
i.e. there exists a morphism 
$$\bmatrix \sigma\\ \delta\endbmatrix :X\to B^\prime\oplus A,$$
such that 
$$ \bmatrix \sigma\\ \delta\endbmatrix \circ (-\beta, f)\ =\ 1_X.$$
This shows that 
$[f]\circ [\delta]=[g]$ in category $\E(\A)$.
and so $(0\to X)$ is projective.

Conversely, suppose that $(\alpha:A^\prime\to A)$ represents a projective
object of $\E(\A)$. We may assume that $\alpha$ is an injective as a map
(by 2.4). The diagram
$$
\CD
(0@>>>A)\\
 & &                @VV{\id}V\\
(A^\prime@>\alpha>>A)
\endCD\tag8
$$
represents an epimorphism in $\E(\A)$ and thus there exists a lifting.
It is represented by a diagram
$$
\CD
(0@>>>A)\\
 & &                @AAfA\\
(A^\prime @>>{\alpha}>A),
\endCD\tag9
$$
where $f:A\to A$ is a morphism in $\H(\A)$ such that 
$$f\circ \alpha\ =0\qquad\text{and}\qquad
f\ =\ 1_A - \alpha\circ g$$ 
for some $g:A\to A^\prime$. It follows that the kernel of $f$ coincides
with the $\im(\alpha)$. Thus we may make an excision with respect to
$P=A^\prime$ and so the original object $(\alpha:A^\prime\to A)$ isomorphic 
to $(0\to A/\alpha(A^\prime))$. This completes the proof.
 
Note additionally, that given an object $(\alpha:A^\prime\to A)$ of $\E(\A)$, 
its projective resolution can be constructed as follows:
$$0\to (0\to A^\prime)@>{[\alpha]}>>(0\to A)@>{[\id]}>>
(\alpha:A^\prime\to A)\to 0$$
(assuming that $\alpha$ is an injective map, cf. 2.4). Thus we conclude that
{\it any object
of $\E(\A)$ admits a projective resolution of length two}. 
Hence {\it the homological dimension of $\E(\A)$ equals to one.}

\subheading{2.9. Example} Consider a chain complex in the abelian category
$\E(\A)$ 
$$\dots\to C_{i+1}@>{\partial}>>C_i @>{\partial}>>C_{i-1}\to \dots,$$
consisting of {\it projective Hilbertian modules} (note that it can be 
equivalently considered as a chain complex in the additive category 
$\H(\A)$). Its $i$-dimensional homology in the extended
abelian category $\E(\A)$ is represented by the morphism
$$H_i(C_\ast)\ =\ (\partial:C_{i+1}\to Z_i)$$
where $Z_i$ is the submodule of cycles, 
$Z_i=\ker [\partial:C_i\to C_{i-1}]$. 

It is an easy exercise, based on Proposition 2.3.

\heading{\bf \S 3. Torsion subcategory $\T(\A)$}\endheading

Now we are going to describe another full subcategory $\T(\A)$ of $\E(\A)$
which is in some sense complementary to $\H(\A)$. We will call it {\it the 
torsion subcategory} and denote $\T(\A)$; its objects will be called
{\it torsion virtual Hilbertian modules}.

\subheading{3.1. Definition} A virtual Hilbertian module 
$(\alpha:A^\prime\to A)$ of $\E(\A)$
will be called {\it torsion} if $A\ =\ \cl(\alpha(A^\prime))$. 

Equivalently,
a virtual Hilbertian module $\X$ is torsion if and only if it admits no 
non-trivial morphisms $\X\to P$ in projective objects $P$ of $\E(\A)$.
The definition in this form clearly does not depend on a particular 
representation of the virtual Hilbertian module and depends only on its
isomorphism class in $\E(\A)$.

Another equivalent definition says: a virtual Hilbertian module represented
by $(\alpha:A^\prime\to A)$ is torsion iff 
$$\dim_\tau(A)\ = \ \dim_\tau(A^\prime)+\dim_\tau(\ker \alpha).$$

$\T(\A)$ will denote the full subcategory of $\E(\A)$
generated by torsion virtual Hilbertian modules.

We will first mention some formal properties of the torsion subcategory.

\proclaim{3.2. Proposition} Given an exact sequence
$$0\to \X^\prime\to \X\to \X^{\prime\prime}\to 0$$
of objects and morphisms of category $\E(\A)$, the middle object $\X$
is torsion if and only if both $\X^\prime$ and $\X^{\prime\prime}$ are
torsion.
\endproclaim
\demo{Proof} Suppose first that both $\X^\prime$ and $\X^{\prime\prime}$
are torsion. Then any morphism $\X\to P$, where $P$ is projective,
vanishes on $\X^\prime$ and so it can be factorized through a morphism
$\X^{\prime\prime}\to P$ which also must vanish, since $\X^{\prime\prime}$
is torsion.

Assume now that $\X$ is torsion. Then clearly $\X^{\prime\prime}$ is torsion
and we are left to show that a {\it subobject of a torsion module is
torsion}. Suppose that 
$$[f]: (\alpha: A^\prime\to A)\to (\beta:B^\prime\to B)\tag10$$
represents a monomorphism in $\E(\A)$ with torsion object
$(\beta: B^\prime\to B)$. As was shown in 2.4 and in 2.6, we may assume 
that all morphisms $\alpha,\beta, f$ are injective. Also,
$\alpha(A^\prime)\supset f^{-1}(\beta(B^\prime))$, cf. 2.5.
Therefore we obtain
$$\dim_\tau(A^\prime)\ge \dim_\tau(B^\prime)\ =\ \dim_\tau(B)\ge 
\dim_\tau(A).$$
Thus the subobject $(\alpha:A^\prime\to A)$ is also torsion.
$\square$
\enddemo

\subheading{3.3} We will see now that any object of the extended category
$\E(\A)$ determines canonically a pair of objects, a torsion and a projective. 

Let $\X=(\alpha: A^\prime\to A)$ be an object of $\E(\A)$. Denote 
$$T(\X)= (\alpha: A^\prime \to \cl(\im(\alpha)))\qquad\text{and}\quad
P(\X)=(0\to A/\cl(\im(\alpha))).$$
$T(\X)$ is clearly a subobject of $\X$ (by 2.7.1), and $P(\X)$ is a 
factor-object of $\X$ (by 2.7.2). It is easy to see that we have a short 
exact sequence
$$0\to T(\X)\to \X\to P(\X)\to 0\tag11$$
in $\E(\A)$. We will say that $T(\X)$ is the {\it torsion part of $\X$} and 
that $P(\X)$ is {\it the projective part of $\X$}. Note that the exact 
sequence (11) splits (since $P(\X)$ is projective) but the splitting is
not canonical.

\subheading{3.4} Given a morphism $f:\X\to \Y$ in the category $\E(\A)$, 
consider the diagram 
$$
\CD
0@>>> T(\X)@>i>>     \X@>j>> P(\X)@>>> 0\\
@.     @V{f^\prime}VV @VfVV   @VV{f^{\prime\prime}}V\\
0@>>> T(\Y)@>{i^\prime}>> \Y@>{j^\prime}>> P(\Y)@>>> 0.$$
\endCD
$$
The morphism $j^\prime\circ f\circ i$ is zero and so there exists a 
unique morphism
$f^\prime : T(\X)\to T(\Y)$ such that $f\circ i=i^\prime\circ f^\prime$.
Thus, we obtain, that {\it any morphism of $\E(\A)$ maps the torsion part
of $\X$ into the torsion part of $\Y$.}

Similarly, morphism $f$ above {\it uniquely determines a morphism
$f^{\prime\prime}:P(\X)\to P(\Y)$ between the projective parts}.

We conclude that {\it there are defined two covariant functors
$$\T: \E(\A)\to \T(\A)\quad\text{and}\quad P:\E(\A) \to \H(\A),$$
which we will call the torsion and the projective part, respectively.}

\subheading{3.5} Given an exact sequence
$$0\to \X^\prime\to \X\to \X^{\prime\prime}\to 0\tag12$$      
in $\E(\A)$, it determines the exact sequence
$$0\to T(\X^\prime)\to T(\X)\to T(\X^{\prime\prime})\tag13$$      
of the torsion parts.
Thus, {\it the functor of torsion part is left exact}. 

Note, that the similar sequence 
$$0\to P(\X^\prime)\to P(\X)\to P(\X^{\prime\prime})\to 0\tag14$$      
for projective parts may be not exact in the middle term.
But this sequence is always {\it weakly exact}, i.e. it is exact if 
considered as a sequence in the original additive
category $\H(\A)$.

It turns out that the homology of the sequences (13) and (14) coincide.
We will formulate this as the following proposition.

\proclaim{3.6. Proposition} Suppose that 
$$0\to \X^\prime\to \X\to \X^{\prime\prime}\to 0$$     
is an exact sequence in $\E(\A)$ and let $\H$ be defined by the exact
sequence
$$0\to T(\X^\prime)\to T(\X)\to T(\X^{\prime\prime})\to \H\to 0$$
Then the homology of the complex (14) of the projective parts in the middle
term is isomorphic to $\H$. \endproclaim
\demo{Proof} Consider the original exact sequence as a chain complex 
$C_\ast$ and let $T_\ast$ denotes its subcomplex formed by the torsion parts.
Similarly, let $P_\ast$ denote the factor-complex formed by the projective
parts. Then we have an exact sequence of chain complexes in $\E(\A)$
$$0\to T_\ast\to C_\ast\to P_\ast\to 0.$$
From the long exact sequence for homology we obtain the isomorphism
$H_i(P_\ast)\to H_{i-1}(T_\ast)$ and in particular 
$H_1(P_\ast)\simeq H_0(T_\ast)=\H$.  $\square$\enddemo

\subheading{3.7. Example} Consider a chain complex in $\E(\A)$ 
$$\dots\to C_{i+1}@>{\partial}>>C_i @>{\partial}>>C_{i-1}\to \dots\tag15$$
consisting of {\it projective Hilbertian modules}.
As we have seen in Example 2.9, its 
$i$-dimensional homology in the extended
abelian category $\E(\A)$ is represented by the morphism
$$H_i(C_\ast)\ =\ (\partial:C_{i+1}\to Z_i),\tag16$$
where $Z_i$ is the submodule of cycles. The projective part of the homology
is the following Hilbertian module 
$$P(H_i(C_\ast))\ =\ Z_i/\cl(\im([\partial:C_{i+1}\to C_i])),\tag17$$ 
which coincides with the definition of reduced $L^2$ homology of $C_\ast$.
Observe, that the full homology $H_i(C_\ast)$, as an object of $\E(\A)$,
may have a non-trivial torsion part
$$T(H_i(C_\ast))\ =\ (\partial: C_{i+1}\to
\cl(\im([\partial:C_{i+1}\to C_i]))).\tag18$$ 
Note also that $T(H_i(C_\ast))$ coincides with the torsion of the virtual
Hilbertian module $(\partial:C_{i+1}\to C_i)$:
$$T(H_i(C_\ast))\ =\ T(\partial:C_{i+1}\to C_i)\tag19$$

We will see in the sequel that the torsion part of the homology of
$C_\ast$ determines 
completely the Novikov-Shubin invariants.

\subheading{3.8. Duality for torsion objects} It turns out that there are two
different notions of duality in the extended abelian category $\E(\A)$, 
one duality for 
projective objects and another duality for torsion objects. The 
duality for projective objects in $\E(\A)$ coincides essentially with the 
duality in
$\H(\A)$, described in 1.6. Now we will construct duality for torsion objects.

Let $\X=(\alpha: A^\prime\to A)$ be a torsion object represented by an 
injective morphism
$\alpha$. Define the {\it dual torsion Hilbertian module}
$\e(\X)$ by
$$\e(\X)\ =\ (\alpha^\ast: A^\ast\to (A^\prime)^\ast),$$
where all $\alpha^\ast$, $A^\ast$ and $(A^\prime)^\ast$ are defined as 
explained in subsection 1.6.

Suppose now that we have two torsion objects $\X=(\alpha: A^\prime\to A)$
and $\Y=(\beta: B^\prime\to B)$ with injective $\alpha$ and $\beta$ and let
$[f]:\X\to\Y$ be a morphism represented by a diagram
$$
\CD
(A^\prime@>\alpha>>A)\\
 & &                @VVfV\\
(B^\prime@>\beta>>B).
\endCD\tag20
$$
According to definition 2.2, there exists a morphism $h:A^\prime\to B^\prime$
making this diagram commutative; this $h$ is in fact unique, because of
injectivity of $\beta$. We define the {\it dual of 
$[f]$} as the morphism
$$\e([f])\ =\ [h^\ast] : \e(\Y)\ \to \e(\X).\tag21$$
It is represented by the diagram
$$
\CD
(B^\ast@>{\beta^\ast}>>B^{\prime\ast})\\
& &      @VV{h^\ast}V\\
(A^\ast@>>{\alpha^\ast}>A^{\prime\ast}).
\endCD
$$
We have to check correctness of this definition. If $F:A\to B^\prime$
is an arbitrary morphism, then the morphism $f^\prime=f+\beta\circ F$
represents the same morphism $[f]=[f^\prime]$ in $\E(\A)$. 
Then the corresponding 
morphism $h^\prime$ is $h^\prime= h+F\circ\alpha$ and thus 
$$h^{\prime\ast}\ =\ h^\ast\ +\ \alpha^\ast\circ F^\ast,$$
which means that $h^\ast$ and $h^{\prime\ast}$ represent the same morphism in
$\E(\A)$.

Note that for any torsion object $\X$ the dual torsion object $\e(\X)$
is isomorphic to $\X$, but not canonically. This follows from the existence
of self-adjoint representation, cf. 4.1.

Clearly,
$$\e(\e(\X))\ \simeq\ \X,\tag22$$
and this isomorphism is canonical.

\proclaim{3.9. Theorem (Universal coefficients)} Let $C_\ast$ be a projective 
chain complex in $\E(\A)$
$$C_\ast:\quad\dots\to C_{i+1}@>\partial>>C_i@>\partial>>C_{i-1}\to\dots
,\tag23$$
and let $C^\ast$ denote the dual cochain complex
$$C^\ast:\quad\dots\to C^\ast_{i-1}@>{\partial^\ast}>>
C^\ast_i@>{\partial^\ast}>>
C^\ast_{i+1}\to\dots,\tag24$$
where the duality is understood as in 1.6.
Consider homology of chain complex $C_\ast$ and cohomology 
of cochain complex $C^\ast$
in the extended abelian category $\E(\A)$. Then the projective part of the
$i$-dimensional cohomology $\H^i(C^\ast)$ is canonically
dual (in the sense of 1.6) to the projective part
of the $i$-dimensional homology $\H_i(C)$, i.e.
$$P(\H^i(C^\ast))\ \simeq \ (P(\H_i(C)))^\ast.\tag25$$
Moreover, the torsion part of the $i$-dimensional cohomology $\H^i(C^\ast)$
is canonically dual (in the sense of 3.8) to torsion part of the 
$(i-1)$-dimensional homology $\H_i(C)$, i.e.
$$T(\H^i(C^\ast))\ \simeq\ \e(T(\H_{i-1}(C))).\tag26$$
\endproclaim
\demo{Proof} As usual, let's denote the space of boundaries by 
$B_i=\im [d:C_{i+1}\to C_i]$, and the space of cycles by
$Z_i=\ker [d:C_{i}\to C_{i-1}]$. Let $H_i$ denote a complement to 
$\overline B_i$ inside $Z_i$ and let $X_i$ be a complement to $Z_i$
inside $C_i$. Thus, we have the decompositions
$$C_i\ =\ \overline B_i \oplus H_i\oplus X_i,\qquad Z_i\ =\ \overline B_i
\oplus H_i$$
and the differential $\partial$ restricts to an injective morphism
$\alpha: X_i\to \overline B_{i-1}$ with dense image. The extended $L^2$
homology of $C$ is represented (cf. 3.7) by
$$
\aligned
\H_i(C)\ =\ & (\partial:C_{i+1}\to Z_i)\ \simeq \\
& (\partial:X_{i+1}\to Z_i)\ \simeq\  \\
& (\alpha: X_{i+1}\to \overline B_i)\ \oplus H_i
\endaligned
$$
Thus, the projective part of $\H_i(C)$ is $H_i$ and the torsion part is
$(\alpha: X_{i+1}\to \overline B_i)$.

Consider now the cochain complex $C^\ast$. We have
$$C_i^\ast\ =\ \overline B_i^\ast\oplus H^\ast\oplus X_i^\ast$$
and $\partial^\ast:C_i^\ast\to C^\ast_{i+1}$ vanishes on 
$H_i^\ast\oplus X_i^\ast$
and it is injective on $\overline B_i^\ast$. Thus, we obtain
$$
\aligned
& \H^i(C^\ast)\ =\ 
(\partial^\ast:C_{i-1}^\ast\to (\ker:d^\ast:C_i^\ast\to C_{i+1}^\ast))\ \simeq\\
& (\partial^\ast:\overline B_{i-1}^\ast \to X_i^\ast\oplus H_i^\ast)\ \simeq\ 
(\alpha^\ast:\overline B_{i-1}^\ast\to X_i^\ast)\oplus H_i^\ast.
\endaligned
$$
Hence, $P(\H^i(C^\ast))=H^\ast_i$ and $T(\H^i(C^\ast))= (\alpha^\ast:
\overline B_{i-1}^\ast\to X_i^\ast)$, which is equal to $\e(T(\H_{i-1}(C)))$,
according to the definition above (cf. 3.8). $\square$ \enddemo

\heading{\bf \S 4. Density functions and the Novikov-Shubin invariants}
\endheading

Von Neumann dimension is a natural numerical invariant of projective
virtual Hilbertian modules. In this section we discuss a numerical invariant
of torsion Hilbertian modules, known as the Novikov-Shu\-bin invariant.
It was introduced by S.P.Novikov and M.A.Shu\-bin in \cite{NS} and then 
studied by M.Gromov and M.A.Shu\-bin \cite{GS}, and also by J.Lott and 
W.L\"uck \cite{LL}; it was considered as an invariant of Hilbert chain
complexes. Our point view in this section is slightly different; we consider
the Novikov-Shubin number (or, more generally, the equivalence class of 
the spectral density function) as an invariant of torsion virtual Hilbertian
modules.

We will see later (cf. \S 5) that the Novikov-Shubin invariant do not 
determine
the isomorphism type of torsion Hilbertian module. We will also construct
some other numerical invariant in the sequel (cf. \S 7).

\subheading{4.1} First we observe that the polar decomposition theorem
(cf., for example, \cite{Di}, appendix 3) implies that 
{\it any torsion virtual Hilbertian 
module $\X$ over $\A$ admits a self-adjoint representation}, i.e. 
representation of the form
$$\X\ =\ (\alpha: A\to A),\tag27$$
where $A$ is a finitely generated Hilbert module over $\A$ (with a fixed
admissible scalar product) and $\alpha$ is a self-adjoint positive operator
($\alpha^\ast=\alpha$, $\alpha >0$), commuting with the action of the von
Neumann algebra $\A$.

Since we will use this fact a few times in this paper, we will present a 
complete proof.

Let $\X=(\beta:A^\prime\to A)$ be an arbitary representation of a torsion
virtual Hilbertian module $\X$ with $\beta$ injective. Then the image
of $\beta$ is dense in $A$. Consider arbitrary admissible scalar products
$\<\ ,\ \>_A$ and $\<\ ,\ \>_{A^\prime}$ on $A$ and $A^\prime$ 
correspondingly (cf. 1.2). This choice determines uniquely a continuous 
self-adjoint positive operator $S:A^\prime\to A^\prime$ such that
$$\<\beta(x),\beta(y)\>_A\ =\ \<Sx,y\>_{A^\prime}$$
for all $x,y\in A^\prime$. It follows that this operator $S$ is injective
and commutes with the action of the von Neumann algebra $\A$; the image of
$S$ is dense.

Let $T\ =\ S^{1/2}: A^\prime\to A^\prime$ be the positive square 
root of $S$. Then $T$ is injective with dense image and commuting with the 
action of $\A$. We obtain that for all $x,y\in A^\prime$
$$\<\beta(x),\beta(y)\>_A\ =\ \<Tx,Ty\>_{A^\prime},$$  
and thus the map $U=\beta\circ T^{-1}$ is an isometry 
$T(A^\prime)\to \im(\beta)$; this isometry can be uniquely extended to
an isometry $U:A^\prime\to A$. 
We obtain that the torsion virtual Hilbertian module 
$$(\alpha=\beta\circ U^{-1}:A\to A)$$
is isomorphic to $(\beta:A^\prime\to A)$ and thus, it gives a self-adjoint 
representation of $\X$.

The above argument also proves the following useful fact:

\proclaim{4.2. Corollary} If $(\beta:A^\prime\to A)$ is a representation of 
a torsion virtual Hilbertian module by an injective morphism $\beta$ then
the Hilbertian modules $A$ and $A^\prime$ are isomorphic.\endproclaim

We may apply the spectral theorem to (27) to get the representation
$$\alpha \ =\ \int_0^\infty \lambda dE_\lambda\tag28$$
where $E_\lambda$ are the self-adjoint projectors which commute with 
$\alpha$ and with the action of the algebra $\A$.

\proclaim{4.3} Given a torsion Hilbertian module $\X$, it admits a 
self-adjoint
representations of the form $(\alpha: A\to A)$ with $\dim_\tau(A)$ being
arbitrarily small.\endproclaim

\demo{Proof} Start from an arbitrary self-adjoint representation
$\X=(\alpha: A\to A)$. For any $\epsilon >0$ there exists $\delta >0$ such that
$\dim_\tau (E_\delta A)<\epsilon$. Denote $A_\delta=E_\delta A$ and let 
$\alpha_\delta :A_\delta\to A_\delta$ acts as $\alpha$.
Then we have 
$$(\alpha: A\to A) \ \simeq (\alpha_\delta: A_\delta\to A_\delta)$$
are isomorphic as objects of $\E(\A)$ (since one is obtained from the other
by excision (cf. 2.4) with respect to the submodule $(1-E_\delta)A\subset A$.
$\square$ \enddemo

\subheading{4.4} Suppose that a torsion Hilbertian module $\X$ 
in a self-adjoint representation $\X=(\alpha: A\to A)$ is given. Denote 
$$F(\lambda)\ =\ \dim_\tau(E_\lambda A),$$
where $E_\lambda$ is the spectral projector determined by (28). Clearly,
$F(\lambda)$ is a monotone non-decreasing right continuous function, defined
for $\lambda\ge 0$; note that $F(\lambda)>0$ for $\lambda>0$ and $F(0)=0$. 
$F(\lambda)$ will be called {\it the
spectral density function of the torsion module $\X$}. This notion 
was introduced by S.P.Novikov and M.A.Shubin \cite{NS}. 

S.P.Novikov and M.A.Shubin \cite{NS} also found an appropriate equivalence
relation between the spectral density functions, called {\it dilatational 
equivalence}, cf. also \cite{GS}, and \cite{LL}. It turns out that up to
this equivalence relation the spectral density function {\it is an invariant
of a torsion Hilbertian module}. 

Let's recall this notion. 
Two spectral density functions $F(\lambda)$ and $G(\lambda)$ are called
{\it dilatationally equivalent} if there exist constants
$C>1$ and $\epsilon>0$ such that 
$$G(C^{-1}\lambda)\le F(\lambda)\le G(C\lambda)$$
holds for all $\lambda\in [0,\epsilon)$.     

\proclaim{4.5. Proposition} The spectral density function, considered up
to dilatational equivalence, is an invariant of isomorphism class
in $\E(\A)$ of a torsion Hilbertian module.\endproclaim

\demo{Proof} First observe that the equivalence class of spectral density 
function $F(\lambda)$
of torsion Hilbertian module $\X=(\alpha:A^\prime\to A)$ does not depend
on the choice of admissible scalar products on $A$ and on $A^\prime$ which
are used in the definition of $F(\lambda)$. This easy follows from results
of \S 1 of \cite{LL}.

Now, considering only injective representations of torsion modules 
(which does not restrict generality by 2.4), we observe that morphisms 
between $(\alpha:A^\prime\to A)$ and $(\beta:B^\prime\to B)$ in the abelian
category $\E(\A)$ (according to definition 2.2) coincide with homotopy 
classes 
of morphisms between $(\alpha:A^\prime\to A)$ and $(\beta:B^\prime\to B)$
(which are now considered as {\it short chain complexes}). Thus, isomorphism 
in $\E(\A)$ corresponds to homotopy equivalence of the chain complexes.
Hence to complete the proof we may use the result of Gromov and Shubin 
\cite{GS}, corollary 2.6, which states equivalence of the spectral density 
functions of chain homotopy equivalent chain complexes.
$\square$
\enddemo

\subheading{4.6. Definition} The Novikov-Shubin invariant of a non-trivial
torsion Hilbertian module $\X$ is defined as
$$\ns(\X)\ =\ \lim_{\lambda\to 0^+}\inf \frac{\ln F(\lambda)}{\ln \lambda}\  
\in [0,\infty]\tag29$$

\subheading{4.7} Observe, that if $\X=\X_1\oplus\X_2$ is direct sum of two
torsion objects, then the spectral density function $F(\lambda)$ of $\X$ is
the sum
$$F(\lambda)\ =\ F_1(\lambda)+F_2(\lambda)\tag30$$
of the spectral density functions $F_1(\lambda)$ and $F_2(\lambda)$ of $\X_1$
and $\X_2$, respectively. It follows that
$$\ns(\X_1\oplus\X_2)\ =\ \min\{\ns(\X_1), \ns(X_2)\}.\tag31$$

\subheading{4.8} The results of \cite{LL} seem to suggest that a more 
convenient numerical invariant of torsion Hilbertian modules $\X$ is
given by
$$\c(\X)\ =\ \ns(\X)^{-1}.$$
I will call this number {\it the capacity of $\X$}. 

Then formula (31) can be 
rewritten as 
$$\c(\X_1\oplus\X_2)\ = \max\{\c(\X_1), \c(\X_2)\}\tag32$$

\proclaim{4.9. Proposition} For a short exact sequence 
of torsion virtual Hilbertian modules 
$$0\to \X_1\to \X\to\X_2\to 0$$
holds
$$\max\{\c(\X_1),\c(\X_2)\}\le\c(\X)\le \c(\X_1)+\c(\X_2).\tag33$$
\endproclaim
\demo{Proof} Using Corollary 2.7, we may represent the inclusion $\X_1\to\X$ 
by a commutative diagram
$$
\CD
(A^\prime@>\alpha>>A)\\
 @VhVV   @VVfV\\
(B^\prime@>\beta>>B)
\endCD
$$
with $h$ isomorphism and $\alpha$, $\beta$, $f$ monomorphisms with dense
images. By Corollary 4.2 all Hilbertian modules $A^\prime$, $A$, $B^\prime$,
$B$ are isomorphic. Thus, we may represent the inclusion $\X_1\to \X$ by
a commutative diagram of the form
$$
\CD
(C@>\alpha>>C)=\X_1\\
 @V{\id}VV   @V{\gamma}VV\\
(C^\prime@>\beta>>C)=\X
\endCD
$$
where $\alpha$, $\beta$, $\gamma$ are injective morphisms $C\to C$ with 
dense images.
Here $\beta=\gamma\circ\alpha$. Applying Proposition 2.3, we obtain that
the cokernel of the map $\X_1\to \X$ can be represented by 
$\X_2=(\gamma:C\to C)$. 

Now our statement follows from Lemma 1.12 of \cite{LL}, which establishes
inequalities between the Novikov-Shubin numbers of the maps $\alpha$,
$\gamma$ and $\beta=\gamma\circ\alpha$. $\square$
\enddemo

Note that Theorem 2.3 of \cite{LL} follows from our Propositions 3.6 and
4.9.

\proclaim{4.10. Corollary} Fix a number $\nu\in [0,\infty]$ and let
$\T_\nu(\A)$ denote the full subcategory of $\E(\A)$ whose class of
objects constitute all
torsion Hilbertian modules $\X$ over $\A$ with capacity less or equal to
$\nu$. Then $\T_\nu(\A)$ is an abelian category.
\endproclaim

Note that $\T_\infty(\A)\ =\ \T(\A)$ and $\T_0(\A)$ is, generally, not empty.
The categories $\T_\nu(\A)$ form a chain of abelian subcategories.

\heading{\bf \S 5. Some examples}\endheading

Here we will discuss some interesting examples demonstrating the properties
of torsion Hilbertian modules.

\subheading{5.1} Consider the group ring $\C[\Z]$ of the infinite cyclic 
group $\Z$. The group $\Z$, lying in $\C[\Z]$, forms an orthonormal base of a 
scalar product $\<\ ,\ \>$ on $\C[\Z]$; completing this scalar product, 
we obtain
a Hilbert space $\ell^2(\Z)$ with the action of the group ring $\C[\Z]$. 
The von
Neumann algebra $\A={\Cal N}(\Z)$ is the commutant of this action. 
The von Neumann trace on $\A$ is given by 
$$\tau(a)\ =\ \<a\cdot 1,1\>,$$
where $1\in \ell^2(\Z)$ is the unit of the group ring $\C[\Z]$ considered as 
an element of $\ell^2(\Z)$. 

We are going to compute examples in the extended abelian category $\E(\A)$.

Using Fourier decompositions, one may identify $\ell^2(\Z)$ with the space of 
square integrable function on the circle $S^1$. Then the group ring $\C[\Z]$
will be identified with the space of Laurent polynomials acting on 
$L^2(S^1)$ by multiplication. The von Neumann algebra $\A={\Cal N}(\Z)$ 
will be
identified with the space $L^\infty(S^1)$ of essentially bounded complex
functions, acting by multiplication on $L^2(S^1)$. The von Neumann trace 
$\tau$ on $\A$ coincides then with the Lebesgue integral:
$$\tau(f)\ =\ \int_{S^1}f(z)dz\quad\text{for}\quad f\in L^\infty(S^1).$$
Here $z$ denotes the coordinate along the circle, $z=\exp(i\theta)$.

\subheading{5.2} Fix a point $z_0=\exp(i\theta_0)$ on the unit circle and a
positive number $\nu>0$. Consider the function
$$f(z)\ =\ |z-z_0|^\nu\ \in L^\infty(S^1).$$
It determines the linear bounded self-adjoint operator 
$$\alpha: L^2(S^1)\ \to \ L^2(S^1)$$
given by the mulplication on $f$, i.e. $\alpha(\phi)=f\phi$.
Obviously, $\alpha$ commutes with the action of the von Neumann algebra
$\A=L^\infty(S^1)$. This operator is injective, but it is not invertible,
since the function $|z-z_0|^{-\nu}$ is not essentially bounded. Denote by
$\X_{\nu,\theta_0}$ the following torsion Hilbertian module
$$\X_{\nu,\theta_0}\ =\ (\alpha: L^2(S^1)\ \to L^2(S^1)).\tag34$$ 

\subheading{5.3} Computing the spectral density function $F(\lambda)$ of 
$\X_{\nu,\theta_0}$, we obtain (using the definition
in 4.4 and also the remarks in 5.1):
$$
\aligned
F(\lambda)\ =\ & \mu\{z\in S^1; |z-z_0|^\nu <\lambda\}\\
=\ & \mu\{z\in S^1; |z-z_0|^2 <\lambda^{2/\nu}\}\\
=\ & \pi^{-1}\cos ^{-1}(1-1/2\lambda^{2/\nu}).
\endaligned\tag35
$$
Here $\mu$ denotes the Lebesgue measure. Using the fact that
$$\lim_{x\to 0}\frac{\cos^{-1}(1-x)}{\sqrt{2x}}\ =\ 1,$$ 
we obtain that the Novikov-Shubin invariant of $\X_{\nu,\theta_0}$ equals
$$\ns(\X_{\nu,\theta_0})\ =\ \frac{1}{\nu},\tag36$$
and for the capacity of $\X_{\nu,\theta_0}$ we obtain
$$\c(\X_{\nu,\theta_0})\ =\ \nu.\tag37$$

\subheading{5.4} Observe that the spectral density function of 
$\X_{\nu,\theta_0}$ depends only on $\nu$ and does not depend on the angle
$\theta_0$. 

Let us prove that {\it for different angles $\theta_0\ne \theta_1$
the torsion Hilbertian modules $\X_{\nu,\theta_0}$ and $\X_{\nu,\theta_1}$
are not isomorphic in $\E(\A)$}. 

In fact, we will now show that for $\theta_0\ne\theta_1$ any morphism
$$\X_{\nu,\theta_0}\ \to \ \X_{\nu,\theta_1}$$
is zero.

It is clear that any morphism $\X_{\nu,\theta_0}\ \to \ \X_{\nu,\theta_1}$
can be presented by a commutative diagram
$$
\CD
L^2(S^1)@>{|z-z_0|^\nu}>> L^2(S^1)\\
@V{m_g}VV                 @VV{m_f}V\\
L^2(S^1)@>>{|z-z_1|^\nu}> L^2(S^1)
\endCD\tag38
$$
where $z_1=\exp(i\theta_1)$ and $m_f$ and $m_g$ are operators of multiplication
by some functions $f,g\in L^\infty(S^1)$. Commutativity of this diagram
means that
$$|z-z_0|^\nu f(z)\ =\ g(z)|z-z_1|^\nu,$$
which can be rewritten as 
$$\frac{g(z)}{|z-z_0|^\nu}\ =\ \frac{f(z)}{|z-z_1|^\nu} \ \equiv h(z).\tag39$$
Now we see that the function $h(z)$ is essentially bounded, 
$h\in L^\infty(S^1)$,
since from the representation as the first fraction in (39) we obtain that 
$h(z)$ is essentially bounded everywhere except a neighbourhood of $z_0$  and,
similarly, from representation as the second fraction in (39) we obtain that
$h(z)$ is essentially bounded everywhere except a neighbourhood of $z_1$.

The multiplication by $h$ defines a morphism $m_h:L^2(S^1)\to L^2(S^1)$ and
$m_f\ =\ m_h \cdot |z-z_1|^\nu$. Now, by the definition of morphisms of
the category $\E(\A)$ (cf. 2.2) we get that the morphism 
$\X_{\nu,\theta_0}\ \to \ \X_{\nu,\theta_1}$
under consideration, vanishes.

\proclaim{5.5. Corollary} The spectral density function and the 
Novikov-Shubin invariant $\ns(\X)$ do not determine the isomorphism type
of a torsion Hilbertian module $\X$.\endproclaim

The constructed example suggests the notion of {\it "local Novikov-Shubin
invariants"} which would measure the spectral density of the $\theta$-{\it
"local part"} of a torsion Hilbertian module for any angle $\theta$. 
We are going to discuss this
subject in another place.

\subheading{5.6} Observe also that for equal angle $\theta_0$ and 
different values of
$\nu$ the torsion modules $\X_{\nu,\theta_0}$ do admit nontrivial maps.
For example, let's mention that for $\nu^\prime>\nu$ there is an exact
sequence of torsion modules
$$0\to \X_{\nu,\theta_0} \to \X_{\nu^\prime,\theta_0}\to 
\X_{\nu-\nu^\prime,\theta_0}\to 0$$
and this sequence does not split (since we have
$\c(\X_{\nu^\prime,\theta_0})
=\c(\X_{\nu,\theta_0})+\c(\X_{\nu-\nu^\prime,\theta_0}),$
cf. 4.7).
         
This gives an example of a torsion module containing a chain of submodules
having continuum members.

\heading{\bf \S 6. Extended $L^2$ homology and cohomology of 
cell complexes}\endheading

Here we will define cohomological functors from the category of polyhedra
to the extended abelian category $\E(\A)$; we will call them {\it extended
$L^2$ homology and cohomology}. We will see that they are homotopy invariant.
The projective part of the extended $L^2$
cohomology coincides with the reduced $L^2$ cohomology; the 
torsion part of the extended $L^2$ cohomology determines the Novikov-Shubin
invariants. Thus, this general homotopy invariant construction gives a very 
transparent (and short!) proof of homotopy invariance
of both, the reduced $L^2$-cohomology (the fact,
proven by J.Dodzuik in \cite{D}), 
and also homotopy invariance of the Novikov-Shubin invariants (established
by M.Gromov and M.A.Shubin in \cite{GS}). 

Note, that the isomorphism type
of the extended $L^2$ cohomology {\it is not determined} by the von Neumann
Betti numbers and the Novikov-Shubin invariants, as examples
in \S 5 show.

We prove in this section the universal coefficients theorem expressing 
the extended $L^2$ homology through cohomology and vice versa. We also prove
the version of Poincar\'e duality.

\subheading{6.1} Let $X$ be a connected CW complex having finitely many
cells in each dimension, and let 
$\pi=\pi_1(X)$ be its 
fundamental group. We denote by $C_\ast(\tilde X)$ the chain complex 
of the universal covering $\tilde X$ of $X$ generated by the lifts of cells 
of $X$. It is a complex of free left $\C[\pi]$-modules.

Let $\A$ be a finite von Neumann algebra with a fixed finite, normal, and
faithful trace $\tau$ and let $M$ be a finitely generated 
Hilbertian module over $\A$ (cf. \S 1). Suppose that a representation
$$\rho:\ \pi\ \to \ \B_{\A}(M)$$
is given; here $\B_{\A}(M)$ denotes the commutant of $M$, i.e. the space of 
all linear bounded operators $M\to M$ commuting with the action of $\A$. We 
will say that any such $M$, carrying the above structures, is a 
$(\A-\pi)$-{\it Hilbertian bimodule}.

The chain complex
$$M\otimes_{\pi}C_\ast(\tilde X)$$
can be considered as a projective complex in the extended abelian 
category $\E(\A)$. Denote by
$$\H_i(X,M),\qquad i=0,1,2\dots$$
the homology of this complex; it is an object of $\E(\A)$. We will
call $\H_i(X,M)$ {\it the extended $L^2$ homology of $X$ with coefficients
in $M$.}  

Observe, that as an object of the extended category, $\H_i(X,M)$ is a 
direct sum of its projective and torsion parts. Using the definition
(cf. 3.3), we see that the projective part $P(\H_i(X,M))$ equals to the 
reduced $L^2$ 
homology (defined as the factor-module
of the cycles of the complex $M\otimes_{\pi}C_\ast(\tilde X)$ divided by
the {\it closure } of the submodule of boundaries). The torsion
part of $\H_i(X,M)$ is isomorphic to the torsion part of the following 
virtual Hilbertian module
$$(d_{i+1}:\ M\otimes_\pi C_{i+1}(\tilde X)\ \to 
M\otimes_\pi C_i(\tilde X)),\tag40$$
cf. 3.7.

To define {\it extended $L^2$ cohomology}, we will 
assume that $M$ satisfies the similar conditions as above, with the  
following modifications: the von Neumann algebra $\A$ acts on $M$
from the right and the group $\pi$ acts on $M$ from the left, so that
$M$ is a {\it $(\pi-\A)$-Hilbertian bimodule}. Then we form the complex
$\Hom_\pi(C_\ast(\tilde X), M).$
The von Neumann alegebra $\A$ acts on this complex from the right; let us
transform this right action into left action by using the involution of
$\A$, i.e. by the rule $\lambda\cdot x=x\cdot \lambda^\ast$.
After this has been done, we obtain a projective cochain complex in the 
extended category $\E(\A)$ which we denote 
$\overline{\Hom_\pi(C_\ast(\tilde X), M)}.$ 
Let
$$\H^i(X,M),\qquad i=0,1,2\dots,\tag41$$
denote the cohomology (understood in $\E(\A)$) of this cochain complex; 
these virtual Hilbertian modules will be called 
{\it extended $L^2$ cohomology of $X$ with coefficients in $M$.}

Again, the projective part $P(\H^i(X,M))$ of the extended cohomology
coincides with the reduced $L^2$ cohomology, which uses
the closure of the space of coboundaries in its definition.

There are two natural constructions
leading to extended $L^2$ homology and cohomology.

\subheading{6.2. Example} Given a connected polyhedron $X$, having finitely
many cells in each dimension, with 
$\pi=\pi_1(X)$, consider the 
$L^2$-completion of the group ring $M=\ell^2(\pi)$ with the right action of 
$\pi$ (right regular representation) and with the left action of the
von Neumann algebra ${\Cal N}(\pi)$ of $\pi$. Recall that the von Neumann
algebra ${\Cal N}(\pi)$ consists of all bounded linear maps 
$\ell^2(\pi)\to \ell^2(\pi)$ commuting with the action of $\pi$ from the
right. 
Thus, we may
form the extended $L^2$ homology
$\H_i(X,\ell^2(\pi))$ (viewing $\ell^2(\pi)$ as a $({\Cal N}(\pi)-\pi)$-Hilbertian
bimodule)
and the extended $L^2$ cohomology $\H^i(X,\ell^2(\pi))$ 
(if we will consider $\ell^2(\pi)$ as a $(\pi-{\Cal N}(\pi))$-Hilbertian
bimodule).

\subheading{6.3. Example} A slightly more general construction consists 
in the following. Suppose that $X$ is a connected polyhedron with finitely
many cells in each dimension, supplied with a finite dimensional 
representation $\pi=\pi_1(X)\to \End_\C(V)$. Consider $V$ as a right 
$\C[\pi]$-module via this representation. Then the tensor product
$\ell^2(\pi)\otimes_\C V$ has a natural $({\Cal N}(\pi)-\pi)$-Hilbertian
bimodule structure,
where the right $\pi$-action is the diagonal one. Thus, we may consider the
extended homology of $X$ with coefficients in $M=\ell^2(\pi)\otimes_\C V$.
This is essentially the previous construction, twisted by a flat finite
dimensional bundle $V$ over $X$. Similar constructions were studied by 
I.M.Singer in \cite{S}.

\proclaim{6.4. Theorem} {\bf (Homotopy invariance)} Let $X$ and 
$Y$ be two connected cell complexes having finitely many cells in each 
dimension, and let $f:X\to Y$ be a homotopy equivalence.
Identify the fundamental groups of $X$ and $Y$ via the induced isomorphism
$f_\ast:\pi =\pi_1(X,x_0)\ \to\ \pi_1(Y,f(x_0))$, where $x_0\in X$ is the 
base point. Then for any $(\A-\pi)$-Hilbertian bimodule $M$ (cf. 6.1), 
the homotopy equivalence
$f$ induces natural isomorphisms of virtual Hilbertian modules of extended
$L^2$ homology
$$f_\ast :\ \H_i(X,M)\ \to\ \H_i(Y,M), \quad i=0,1,2,\dots.$$
Similarly, for any $(\pi-\A)$-Hilbertian bimodule $M$, the homotopy 
equivalence
$f$ induces natural isomorphisms of virtual Hilbertian modules of extended
$L^2$ cohomology
$$f^\ast :\ \H^i(Y,M)\ \to\ \H^i(X,M),\quad i=0,1,2,\dots.$$
\endproclaim
\demo{Proof} The induced chain map 
$f_\ast : C_\ast(\tilde X)\to C_\ast(\tilde Y)$ is a chain homotopy 
equivalence. Applying the functor $M\otimes_\pi$, we obtain that the
chain map 
$\id\otimes f_\ast : M\otimes_\pi C_\ast(\tilde X)\to 
M\otimes_\pi C_\ast(\tilde Y)$ is a chain homotopy equivalence between
chain complexes in abelian category $\E(\A)$. Then the induced map
on the homology is an isomorphism. 

Similarly one gets the result concerning
cohomology. $\square$ \enddemo

Isomorphism of Hilbertian modules induces isomorphism 
of their projective and torsion parts. Thus, the above theorem essentially 
states that
both the projective and torsion parts of the extended $L^2$ homology and 
cohomology are homotopy invariant. For the projective part this gives the 
result of J.Dodziuk \cite{D}. For the torsion part this fact combined with
Proposition 4.5
implies homotopy invariance of the Novikov-Shubin numbers and the spectral
density functions (cf. \cite{GS}). 

\proclaim{6.5. Corollary (Homotopy invariance of Novikov-Shubin invariants)}
Let $X$ be a connected polyhedron having finitely many cells in each 
dimension, with fundamental group $\pi$,
and let $M$ be a Hilbertian 
$(\A-\pi)$-bimodule (cf. 6.1). 
Then for any $i=0,1,2,\dots$ the spectral density function (considered up to
dilatational equivalence) and the Novikov-Shubin
number of the torsion submodule of 
$$(\id\otimes d:\ M\otimes_\pi C_{i+1}(\tilde X)\ 
\to\  M\otimes_\pi C_i(\tilde X))\tag42$$
are homotopy invariants of $X$.
\endproclaim

Note, that the original definition of the Novikov-Shubin numbers, given
in \cite{NS}, uses the spectral density functions of the Laplacians, acting
on smooth forms. A.Efremov \cite{E} proved equivalence of the above 
{\it "analytic"} definition to the {\it "combinatorial"} definition which
uses the spectral density functions of the Laplacians acting on the chain
complexes, constructed by the cell decomposition of $X$. J.Lott and 
W.L\"uck \cite{LL} showed that instead of the Laplacian, one may consider
the boundary homomorphism itself.

Comparing our notations with notations of J.Lott and W.L\"uck in \cite{LL},
we note the that there is a shift of dimensions: 
number $\ns(T(\H_i(C)))$ in notations of \cite{LL} is $\alpha_{i+1}(C)$. 

The extended $L^2$ homology and cohomology determine each other as the 
following theorem shows:

\proclaim{6.6. Theorem (Universal coefficients)} Let $X$ be a 
polyhedron with finitely many simplexes in each dimension, 
and let $M$ be a Hilbertian $(\A-\pi)$-bimodule. Then the dual 
module $M^\ast$ is naturally defined as Hilbertian $(\pi-\A)$-bimodule and 
there
are canonical isomorphisms for the projective and torsion parts of the 
extended $L^2$ homology and cohomology
$$P(\H^i(X,M^\ast))\ \simeq \ P(\H_i(X,M))^\ast,\tag43$$
$$T(\H^i(X,M^\ast))\ \simeq \ \e(T(\H_{i-1}(X,M))).\tag44$$
\endproclaim
\demo{Proof} This follows from Theorem 3.9 applied to the projective
chain complex $C\ =\ M\otimes_\pi C_\ast(\tilde X)$. Then the dual cochain
complex is $C^\ast\ =\ \overline{\Hom_\pi(C_\ast(\tilde X),M^\ast)}$
and the statement follows. $\square$\enddemo

There is also Poincar\'e duality relating the extended $L^2$ homology 
and cohomology of complementary dimensions. Here is the precise statement.

\proclaim{6.7. Theorem (Poincar\'e duality)} Let $X$ be an $n$-dimensional
closed PL manifold (or, more generally, Poncar\'e complex) and let 
$w:\pi=\pi_1(X)\to \{+1,-1\}$ denote the first Stifel-Whitney class of
$X$. Let $M$ be a Hilbertian $(\A-\pi)$-bimodule. Denote by $M^w$ the 
following $(\pi-\A)$-Hilbertian bimodule structure on $M$: the right action
of $\pi$ on $M$ is given by
$g\cdot m\ =\ w(g)m\cdot g^{-1}$ where $g\in\pi$, $m\in M,$
and the left action of $\A$ on $M$ is given by
$m\cdot a=a^\ast\cdot m$ for $a\in\A$, $m\in M.$
Then there exists a natural isomorphism of virtual Hilbertian modules
$$\H_i(X,M)\ \to\ \H^{n-i}(X,M^w),\quad i=0,1,2,\dots, n.\tag45$$
\endproclaim
\demo{Proof} Denote by $C$ the chain complex of the universal covering
$\tilde X$ of $X$, $C=C_\ast(\tilde X)$. It is a chain complex of free
finitely generated left $\Z[\pi]$-modules. Let $C^\bullet$ be the following
cochain complex $C^\bullet=\Hom_\pi(C,\Z[\pi])$ of free finitely generated 
right $\Z[\pi]$-modules. Let $(C^\bullet)^D$ denote the chain complex of 
left $\Z[\pi]$-modules, which is obtained from $C^\bullet$, first, by 
enumeration of the dimensions in the opposite way: 
$(C^\bullet)^D_i=(C^\bullet)^{n-i}$ and, secondly, by transforming the right
action of $\pi$ into a left action by using the following involution
of the group ring $\Z[\pi]$: $g\mapsto w(g)g^{-1}$ for $g\in \pi$. 

The Poincar\'e duality theorem (cf. for example, \cite{M}) states that there
exists natural chain homotopy equivalence
$C\ \to \ (C^\bullet)^D$. Tensoring it with $M$, we obtain a chain homotopy
equivalence
$$M\otimes_\pi C\ \to\ M\otimes_\pi(C^\bullet)^D$$
between projective chain complexes in $\E(\A)$. The $i$-dimensional
extended $L^2$ homology of the last chain complex can be identified with 
$\H^{n-i}(X,M^w)$. This implies the theorem. $\square$\enddemo

Combining this theorem with Theorem 6.6 we obtain:

\proclaim{6.8. Corollary} If $X$ as an $n$-dimensional closed PL manifold (or
Poincar\'e complex) with first Stiffel-Whitney class $w$ then for arbitrary
Hilbertian $(\A-\pi)$-bimodule $M$ there are natural isomorphisms for the
projective and torsion parts
$$P(\H_i(X,M))\ \simeq\ P(\H_{n-i}(X,M^{\ast w}))^\ast,\quad 
i=0,1,2,\dots n,\tag46$$
$$T(\H_i(X,M))\ \simeq\ \e(T(\H_{n-i-1}(X,M^{\ast w}))),\quad
i=0,1,2,\dots, n-1.\tag47$$
The first isomorphism uses the duality for projective objects (cf. 1.6),
while the second isomorphism uses the duality for torsion objects (cf. 3.8).
Here $M^{\ast w}$ denotes the Hilbertian $(\A-\pi)$-bimodule which is obtained
from the Hilbertian $(\pi-\A)$-bimodule $M^\ast$ by using the involution of $\A$
and the involution on $\pi$ determined by $w$ (cf. above).
\endproclaim

\subheading{6.9} Note that $T(\H_n(X,M))=0$ for any $n$-dimensional
polyhedron $X$.

\subheading{6.10} Another observation: $\H^0(X,M)=0$ for any Hilbertian
$(\pi-\A)$-bimodule $M$ satisfying the condition: 
there is no $m\in M$, $m\ne 0$, such that $g\cdot m=m$
for any $g\in\pi$. In fact, it is easy to see that $\H^0(X,M)$ has no 
torsion and it is isomorphic to  
$$\{m\in M;gm=m \quad\text{for all}\quad g\in\pi\}$$
which has natural structure of Hilbertian module. It is not difficult to
find examples when the last module is nonzero.

It follows, in particular, that $\H^0(X,\ell^2(\pi))\ =\ 0$ if $\pi$ is 
infinite.

On the contrary, the extended $L^2$ homology in dimension zero 
$\H_0(X,M)$ has non-trivial torsion part in many cases. More precisely,
theorem of R.Brooks \cite{B} states that $\H_0(X,\ell^2(\pi))$ {\it 
has nontrivial torsion part if and only if the fundamental group $\pi$
is amenable}.

\heading{\bf \S 7. Minimal number of generators}\endheading

In this section we will introduce and study 
a new numerical invariant of Hilbertian
modules, which we call the minimal number of generators. 
Considered on projective modules, it is comparable with the von 
Neumann dimension,
although it is always integral and may be much larger, than the dimension.
Unlike the von Neumann dimension, this invariant is also nontrivial on 
torsion modules. We show in this section (by computing examples) that 
this new invariant is independent from the Novikov-Shubin invariant.

We will apply this invariant in the next section to the 
Morse theory.

We denote by $\A$ a fixed finite von Neumann
algebra supplied with a finite, normal and faithful trace $\tau$. We 
consider categories $\H(\A)$ and $\E(\A)$, introduced in \S 1 and \S 2.
Recall, that they actually depend on $\tau$ as well, although in our notation
it is suppressed.

\subheading{7.1} First, observe 
that the Hilbertian 
module $\ell^2(\A)$ is a fixed {\it generator} of the category $\E(\A)$ in the 
sense of the category theory \cite{F}. It means, in this particular case,
that any object $\X$ of category $\E(\A)$ is factor-object of a finite 
direct sum of the form $\oplus \ell^2(\A)$. 

Recall that a Hilbertian module is called {\it free} if it is a direct sum
of several copies of $\ell^2(\A)$. A chain complex is called {\it free} if it
consists of free modules.

We arrive at the following definition.

\subheading{7.2. Definition} Let $\X$ be a virtual Hilbertian module. 
We will denote by $\mu(\X)$ the minimal integer $\mu$ such that there
exists an epimorphism of the direct sum of $\mu$ copies of $\ell^2(\A)$ 
onto $\X$. We will call $\mu(\X)$ {\it the minimal number of generators
of $\X$}.

Obviously, $\mu(\X)=0$ if and only if $\X=0$.

If $\X$ is projective, then $\mu(\X)\ge \dim_\tau (\X)$. 

We also have the following property
$$\max\{\mu(\X),\mu(\Y)\}\ \le \mu(\X\oplus \Y)\ \le\  \mu(\X)\ +\ 
\mu(\Y).\tag48$$
Examples below show that both extreme cases allowed by this inequality can
be realized.

\subheading{7.3. Example} Let $\A={\Cal N}(\Z)$ be the von Neumann algebra of
the infinite cyclic group, cf. 5.1. Let $\X_{\nu, \theta}$ be the torsion
module constructed in 5.2. Then obviously 
$$\mu(\X_{\nu, \theta})\ =\ 1.\tag49$$ 

\subheading{7.4. Example} Here we will present a more interesting example.
Fix an arbitrary $\nu>0$ and an angle $\theta$.
Let $\X$ denote the direct sum of $n$ copies of $\X_{\nu, \theta}$.

\proclaim{Claim} The minimal number of generators of $\X$ equals to $n$.
\endproclaim

\demo{Proof} Let $F$ denote the direct sum of $n$ copies of $\ell^2(\A)$.
Obviously, $\X$ has a representation $\X=(\alpha:F\to F)$, where $\alpha$ is
given by a diagonal $n\times n$-matrix with the functions $|z-z_0|^\nu$
standing along the diagonal (here $z_0=\exp(i\theta)$). We identify
the von Neumann algebra ${\Cal N}(\Z)$ with $L^\infty(S^1)$ and $\ell^2(S^1)$
with $L^2(S^1)$, cf. \S 5. $z$ denotes the coordinate along the circle.

We have to show that there is no epimorphism $[f]:F^\prime \to \X$ if 
$F^\prime$ is a direct sum of $m$ copies of $\ell^2(\A)$ with $m<n$. 
Suppose that such epimorphism exists; it is then represented by a diagram
$$
\CD
  @. F^\prime\\
@.   @VVfV\\ 
\X=(\alpha:F@>>>F).
\endCD
$$
By 2.5, we obtain that the morphism 
$[\alpha,f]:\ F\oplus F^\prime\ \to F$
is an epimorphism in $\H(\A)$. Then there exists a splitting (cf. 1.5)
$$F\oplus F^\prime@<{\bmatrix\beta\\ \gamma\endbmatrix}<< F,$$
so that 
$$\alpha\beta + f\gamma\ =\ \id_F.\tag50$$

Note that the morphisms $\alpha$, $\beta$, $\gamma$, $f$ are represented by
rectangular matrices with entries in the von Neumann algebra 
$\A=L^\infty(S^1)$. Denote by $\chi: L^\infty(S^1)\to \C$ any multiplicative
homomorphism, which has the following property: for any continuous function 
$\phi\in C(S^1)\subset L^\infty(S^1),$ one has $\chi(\phi)=\phi(z_0)$. 
Such multiplicative 
functional exists by virtue of \S 12 of \cite{GRS}.

Now, applying $\chi$ to the matrix equation (50), we obtain 
(since $\chi(\alpha)=0$) the equality
$$\chi(f)\chi(\gamma)\ =\ 1_{n\times n}.\tag51$$
In the last equation, $\chi(f)$ is a $n\times m$ complex matrix and
$\chi(\gamma)$ is an $m\times n$-matrix; (51) is impossible, if $m<n$.
$\square$\enddemo

\subheading{7.5. Example} Suppose now that the angles $\theta_1$ and 
$\theta_2$ are different. Let $\nu$ and $\nu^\prime$ be two arbitrary 
positive numbers. 
\proclaim{Claim}
$$\mu(\X_{\nu,\theta_1}\ \oplus \ \X_{\nu^\prime,\theta_2})\ =\ 1.$$
\endproclaim

To prove this it is enough to show that the diagram
$$
\CD
@. L^2(S^1)\\
@. @VV{\bmatrix \beta\\ \alpha\endbmatrix}V\\
(\bmatrix \alpha 0\\ 0 \beta\endbmatrix: L^2(S^1)\oplus L^2(S^1) 
@>>>L^2(S^1)\oplus L^2(S^1))
\endCD\tag52
$$
represents an epimorphism in category $\E(\A)$.
Here $\alpha$ denotes the operator of multiplication by $|z-z_1|^\nu$,
with $z_1=\exp(i\theta_1)$, and $\beta$ denotes the operator of 
multiplication by $|z-z_2|^{\nu^\prime}$, with $z_2=\exp(i\theta_2)$.
Decompose the circle into the union of two intervals, such that one 
of the intervals contains the 
point $z_1$ in its interior, and the other interval contains $z_2$ in its
interior. Let $\chi_1(z)$ be the characteristic function of the interval
containing $z_1$, and let $\chi_2(z)$ be the characteristic function of the 
other interval. Then we have, $\chi_1(z)+\chi_2(z)=1$ for all $z$ 
(except for two end points of the intervals) and
$\chi_j^2=\chi_j$, where $j=1,2$.

To show that the above diagram represents an epimorphism, we may use the 
criterion 2.5. Thus, we have to show that, given arbitrary pair
$(\phi,\psi)\in L^2(S^1)\oplus L^2(S^1)$, there exist a pair 
$(a,b)\in L^2(S^1)\oplus L^2(S^1)$
and a  function $c\in L^2(S^1)$ such that
$$\bmatrix \alpha 0\\ 0 \beta\endbmatrix \bmatrix a \\ b\endbmatrix \ +\ 
\bmatrix \beta \\ \alpha \endbmatrix \cdot c \ =\ \bmatrix \phi\\ \psi
\endbmatrix$$
This can be achieved by setting
$$
\aligned
a\ &=\ \alpha^{-1}\chi_2\phi \ +\ \beta\alpha^{-2}\chi_2\psi\\
b\ &= \ \beta^{-1}\chi_1\phi\ +\ \alpha\beta^{-2}\chi_1\psi\\
c\ &=\ \beta^{-1}\chi_1\phi \ +\ \alpha^{-1}\chi_2\psi
\endaligned
$$
This proves our statement.

\subheading{7.6} Comparing the numerical invariants $\ns(\X)$ and 
$\mu(\X)$, we conclude {\it that they are independent}. In fact, the sum 
of $n$ copies
of $\X_{\nu,\theta}$ has capacity $\nu$ and its minimal number of generators 
is $n$. Thus, {\it any pair of numbers $\nu\in \R^+$, $n\in \N$ can be 
realized}.

\heading{\bf \S 8. Application: Morse inequalities}\endheading

In the original papers of S.P.Novikov and M.A.Shubin \cite{NS}, \cite{NS1}
(cf. also \cite{Sh})
it was shown how one may use the notion of dimension in the sense of von 
Neumann, in order to improve the classical Morse inequalities for 
the numbers of critical points of functions on compact manifolds.
We are going to show in this section that the phenomenon 
responsible for the Novikov-Shubin invariants (namely, measuring the "rate"
of zero being in the continuous spectrum), which was discovered in the same 
papers \cite{NS}, \cite{NS1} of S.P.Novikov and M.A.Shubin,
may also be used to further strengthening the Morse inequalities. 
More precisely, we will use the numerical invariant $\mu(\X)$,
introduced in the previous section, 
and prove a version of Morse inequalities, which involves this new
invariant applied to the extended $L^2$ homology. 
This approach extracts a quantitative information from the 
torsion part of the extended $L^2$ homology, as well.

\proclaim{8.1. Theorem} Let $\A$ be a finite von Neumann algebra supplied 
with a finite normal and faithful trace $\tau$, and let
$$C:\quad \dots \to C_{i+1}\to C_i\to C_{i-1}\to \dots\tag53$$
be a free chain complex (cf. 7.1) in $\E(\A)$. 
Then for any integer $i$ the following inequality holds:
$$\dim_\tau(C_i)\ \ge \ \mu[\H_i(C)\oplus T(\H_{i-1}(C))].\tag54$$
\endproclaim

\demo{Proof} If $Z_i$ denotes the submodule of $i$-dimensional cycles, then
we can find a splitting
$$C_i\ =\ Z_i\ \oplus \ X_i.$$
Let $f:C_i\to X_i$ and $f^\prime :C_i\to Z_i$ denote the projections.
The restriction of the boundary homomorphism onto $X_i$ determines
an injective bounded linear map $\alpha:X_i\to \overline B_{i-1}$ with
dense image. Thus we obtain
$$T(\H_{i-1}(C))\ =\ (\alpha:X_i\to \overline B_{i-1}).$$
By Corollary 4.2 we get that 
there exists an {\it isomorphisms} of Hilbertian modules 
$g:X_i\to \overline B_{i-1}$. Now, the map
$$\bmatrix f^\prime \\ gf \endbmatrix : C_i\ \to\ Z_i\oplus \overline B_{i-1}$$
is an epimorphism. Composing it with the obvious epimorphism
$$Z_i\oplus \overline B_{i-1} \to \H_i(C)\oplus T(\H_{i-1}(C))=
(C_{i+1}\to Z_i)\oplus (X_i\to \overline B_{i-1}),$$
we obtain an epimorphism from free module $C_i$ onto the last module.
This shows that 
$$\mu[\H_i(C)\oplus T(\H_{i-1}(C))]\ \le \dim_\tau(C_i),$$
and completes the proof. $\square$.
\enddemo

\proclaim{8.2. Theorem (Morse type inequalities)} Let $X$ be a closed
manifold and let $f:X\to \R$ be a non-degenerate Morse function on $X$.
Suppose that a finite dimensional representation $\pi \to \End_\C(V)$
of the fundamental group $\pi=\pi_1(X)$ is given, and let $M$ denote the
Hilbertian $({\Cal N}(\pi)-\pi)$-bimodule $\ell^2(\pi)\otimes_\C V$,
constructed in 6.3. Then for the
Morse numbers $m_i(f)$ of critical points of $f$ of index $i$ 
the following inequalities hold:
$$m_i(f)\ \ge\ (\dim_\C V)^{-1}\cdot 
\mu[\H_i(X,M)\oplus T(\H_{i-1}(X,M))], \quad i=0,1,2,\dots \tag55$$
\endproclaim
\demo{Proof} The theorem follows by applying the previous Theorem 8.1
to the 
chain complex $M\otimes_\pi C_\ast(\tilde X)$, where $C_\ast(\tilde X)$
is constructed by means of the cell decomposition of $X$, determined by the 
Morse function $f$. $\square$
\enddemo

\subheading{8.3. Example} Consider the simplest possible example: 
let $X$ be
the circle $S^1$. In this case all von Neumann Betti numbers vanish and so 
the Morse type inequalities of Novikov and Shubin \cite{NS}, \cite{NS1}
do not predict critical points. 

Let us apply Theorem 8.2 with $V$ the trivial one-dimensional representation. 
Then the only non-vanishing extended homology is
$\H_0(X,M)\ =\ \X_{\nu,\theta}$ (using the notation introduced in 5.2), 
where $\nu\ =\ 1$ and $\theta\ =\ 0$. Applying
Theorem 8.2 we obtain (cf. 7.3) that any Morse function on the circle 
must have at least one minimum and at least one maximum!

\heading{\bf \S 9. Group homology and cohomology 
with values in $\E(\A)$}\endheading

We are going to define in this section some homological functors 
which associate to a representation of a discrete group a sequence of 
objects of the extended 
abelian category $\E(\A)$. The projective part of these functors coincides
with the reduced $L^2$ cohomology, studied, by J.Cheeger
and M.Gromov in \cite{CG}.
We denote these new functors $\TOR^{\pi}_q(N,M)$ and $\EXT_\pi^q(N,M)$ since 
their construction is  
similar to building the usual $\Tor$ and $\Ext$ functors. 

Our main goal in this section is to express the extended $L^2$ homology 
and cohomology of a cell complex $X$ (defined in \S 6), through the homology 
of its universal covering of $X$;
it is done with the aid of $\TOR$ and $\EXT$ functors and a spectral sequence 
in $\E(\A)$.

\subheading{9.1} Let $\A$ be a finite von Neumann algebra supplied with a
finite, normal, and faithful trace $\tau$.

Let $\pi$ be a discrete group.

Let $M$ be a Hilbertian $(\A-\pi)$-bimodule (cf. 6.1),
and let $N$ be a left $\C[\pi]$-module, having a free $\C[\pi]$-resolution
$$\dots\to C_{q+1}\to C_q\to C_{q-1}\to \dots \to C_0 
@>\epsilon>>N\to 0\tag56$$
such that $C_q$ is finitely generated over $\C[\pi]$ in each dimension $q$.
Applying functor $M\otimes_\pi$ to (56), we obtain the following projective
chain complex in $\E(\A)$
$$\dots\to M\otimes_\pi C_{q+1}\to M\otimes_\pi C_q\to M\otimes_\pi 
C_{q-1}\to \dots .\tag57$$
We denote by 
$$\TOR_q^\pi(M,N)\tag58$$
its $q$-dimensional homology in $\E(\A)$. Since any two free resolutions
of $N$ are homotopy equivalent, we obtain that $\TOR_q^\pi(M,N)$ is correctly
defined and is a covariant functor of $N$.

Note, that if $N$ is {\it finitely generated free over} $\C[\pi]$ then the
tensor product $M\otimes_\pi N$
is well defined as a projective object of $\E(\A)$, and 
$\TOR_0^\pi(M,N)\ \simeq \ M\otimes_\pi N$, and $\TOR_p^\pi(M,N)=0$ 
for $p\ge 1$.

\subheading{9.2. Example} Suppose that $\pi=\Z$ and $N=\Z$ with the 
trivial $\pi$-action. Then
$$0\to \C[\Z]@>{z-1}>>\C[\Z]\to \Z\to 0$$
is a free resolution of $\Z$ and thus we obtain
$$\TOR_0^\pi(\ell^2(\Z),\Z)\ =\ \X_{1,0}\tag59$$
(in the notations of 5.2) and
$$\TOR_p^\pi(\ell^2(\Z),\Z)\ =\ 0$$
for all $p\ge 1$.

\subheading{9.3} To define a similar cohomological notion, suppose that
$M$ is a Hilbertian $(\pi-\A)$-bimodule and that $N$ is (as before)
a left $\C[\pi]$-module,
having free resolution (56), which is finitely generated in every dimension.
Form the cochain complex $\Hom_\pi(C_\ast, M)$; it is a cochain complex
of right $\A$-Hilbertian modules. We may convert it into a cochain complex of
left $\A$-Hilbertian modules by using the involution of $\A$; the resulting
complex let's denote by $\overline{\Hom_\pi(C_\ast, M)}$. Now, define
$$\EXT^q_{\pi}(N,M)\ =\ H^q(\overline{\Hom_\pi(C_\ast, M)}).\tag60$$
It is a contravariant functor of $N$.

\subheading{9.4. Example} Suppose that the $\C[\pi]$-module
$N=\Z$ with the trivial $\pi$-action 
admits a free resolution which is finitely generated in every dimension. 
Then the above 
construction produces the {\it group cohomology with values in $\E(\A)$}:
$$\H^p(\pi, M)\ =\ \EXT^p_\pi (M,\Z).\tag61$$
Note that it can be understood also in the framework of construction 
of subsection 6.1 as $\H^p(K(\pi,1),M)$, where $K(\pi,1)$ denotes the 
Eilenberg-MacLane space. In the case, when $M=\ell^2(\pi)$, the projective part
of the cohomology $P(\H^p(\pi, M))$ is denoted in \cite{CG} 
by $\overline H^p_{(2)}(\pi)$.

\subheading{9.5} Using theorem 6.6 we get the following duality relations:
$$P(\EXT^i_\pi(N,M))\ \simeq\ P(\TOR_i^\pi(M^\ast,N))^\ast\tag62$$
and 
$$T(\EXT^i_\pi(N,M))\ \simeq\ \e(T(\TOR_{i-1}^\pi(M^\ast,N))).\tag63$$
Here $M$ is an arbitrary Hilbertian $(\pi-\A)$-bimodule.

\proclaim{9.6. Corollary} The module $\EXT^0_\pi(N,M)$ is always projective
and it is isomorphic to $\overline{\Hom_\pi(N,M)}$.\endproclaim 

The following theorem is one of the main result of this section.
It establishes a relation between the homology of the universal covering,
considered as a module over the group ring of the fundamental group, 
and the extended $L^2$ homology and cohomology.

\proclaim{9.7. Theorem} Let $X$ be a finite cell complex with fundamental
group $\pi=\pi_1(X)$ and let $\tilde X$ denote the universal cover of
$X$. Suppose
that the homology modules $H_q(\tilde X,\C)$, considered as left 
$\C[\pi]$-modules, have free resolutions with finitely generated 
$\C[\pi]$-modules of chains
in all dimensions. Suppose that $M$ is a Hilbertian $(\A-\pi)$-bimodule.
Then there exists a spectral sequence in the abelian 
category $\E(\A)$ with the initial $E^2$-term
$$E^2_{p,q}\ =\ \TOR^\pi_q(M,H_p(\tilde X))\ \Longrightarrow 
\ \H_{p+q}(X,M).\tag64$$
The limit term $E^\infty_{p,q}$ coincides with $E^r_{p,q}$ for some
large $r$.
\endproclaim
\demo{Proof} First, we will construct a special {\it Cartan-Eilenberg 
resolution} $D_{\ast \ast}\to C_\ast$ of the chain complex 
$C_\ast=C_\ast(\tilde X)$, cf. \cite{W}, \S 5.7. We will need this resolution
to satisfy some properties, additional to those, mentioned in \cite{W}.
Here $D_{pq}$ is a double complex, with two differentials
$d^h$ (horizontal) and $d^v$ (vertical). It is required to 
satisfy the following conditions:

$\bullet$ $D_{pq}$ is free and finitely generated over $\C[\pi]$ for any
pair $p,q$;

$\bullet$ $D_{pq}=0$ if $p<0$ or if $p>\dim X$;

$\bullet$ For any number $p$, the modules of boundaries 
$B_p(D_{\ast\ast},d^h)$, 
the module of cycles $Z_p(D_{\ast\ast},d^h)$, and the module
of homology $H_p(D_{\ast\ast},d^h)$
of the double complex with respect to the horizontal differential $d^h$,
are free finitely generated resolutions of $B_p(C)$, $Z_p(C)$,
$H_p(C)$, respectively.

Such Cartan-Eilenberg resolution can be constructed inductively using
the Horseshoe lemma (cf. \cite{W}, 2.2.8) similarly to proof of Lemma 5.7.2
in \cite{W}; in fact we use a version the Horseshoe  
lemma, where instead of projective resolutions we deal with {\it free and 
finitely generated} ones. 

If we are given resolutions of $B_p(C)$ and $H_p(C)$, then, using the 
Horseshoe lemma applied to the exact sequence 
$$0\to B_p(C)\to Z_p(C)\to H_p(C)\to 0,$$
we obtain a resolution of $Z_p(C)$; the resolution of $Z_p(C)$ produces a
resolution of $B_{p-1}(C)$ via the exact sequence 
$$0\to Z_p(C)\to C_p\to B_{p-1}(C)\to 0.\tag65$$
And, at last, using exact sequence (65) again and the constructed 
resolutions of $Z_p(C)$ and of $B_{p-1}(C)$, we obtain a resolution 
of $C_p$. It gives a column
$D_{p,\ast}$ of the Cartan-Eilenberg resulution, standing above $C_p$.
We apply this procedure inductively,
starting from the maximal dimension $p=\dim(X)$ 
and moving downwards.

Suppose now that the Cartan-Eilenberg resolution $D_{\ast\ast}$ with the above
properties has been constructed. 
Consider the following double complex $M\otimes_\pi D_{pq}$ as a double
complex in the abelian category $\E(\A)$. Computing first homology of 
this double complex with respect to vertical differential, we obtain that
that nontrivial homology will appear only on the row $q=0$ and the homology
at point $(p,0)$ is $M\otimes_\pi C_p$. Thus we obtain that {\it the homology
of the total complex $\Tot_\ast(D_{\ast\ast})$ is precisely the extended 
$L^2$ homology $\H_\ast(X,M)$}. 

On the other hand, if we compute first homology of the double complex
$M\otimes_\pi D_{pq}$ with respect to the horizontal differential $d^h$ and
then with respect to the vertical $d^v$, we will get on place $(p,q)$
the following virtual Hilbertian module
$$E^2_{p,q}\ =\ \TOR_q^\pi(M, H_p(\tilde X)).$$
This gives the desired spectral sequence. 

It clearly stabilizes after a finite number of steps.
$\square$.
\enddemo

There is also a cohomological version of Theorem 9.7.

\proclaim{9.8. Theorem} Let $X$ be a finite cell complex with fundamental
group $\pi=\pi_1(X)$ and let $\tilde X$ denote the universal cover. Suppose
that the homology modules $H_p(\tilde X,\C)$ 
(considered as left $\C[\pi]$-modules)    
have free resolutions with finitely generated modules
in all dimensions. Let $M$ be a Hilbertian $(\pi-\A)$-bimodule.
Then there exists a spectral sequence in the abelian 
category $\E(\A)$ with the initial $E_2$-term
$$E_2^{p,q}\ =\ \EXT_\pi^q(H_p(\tilde X),M)\ \Longrightarrow 
\ \H^{p+q}(X,M).\tag66$$
The limit term $E_\infty^{p,q}$ coincides with $E_r^{p,q}$ for some large
$r$.
\endproclaim
\demo{Proof} It is similar to Proof of Theorem 9.7.$\square$\enddemo

\subheading{9.9} As an application, consider the case when the fundamental
group $\pi=\pi_1(X)$ is {\it free}. Note, that the group ring $\C[\pi]$
of the free group is a FIR (cf. \cite{C}) and it is {\it coherent}, 
cf. \cite{C}, page 554. This implies that all $\TOR^\pi_q$ with $q>1$ vanish 
and for any finite polyhedron $X$ with $\pi_1(X)=\pi$ the homology groups 
of the universal covering $H_p(\tilde X)$ are finitely presented as
$\C[\pi]$-modules (so Theorems 9.7 and 9.8 can always be applied).
Applying Theorem 9.7, we obtain the following exact sequence
$$0\to \TOR_1^\pi(M,H_{p-1}(\tilde X))\ \to \ \H_p(X,M)\ \to\ 
\TOR_0^\pi(M,H_p(\tilde X))\ \to 0.\tag67$$
Observe, that the module on the left is projective. 

\proclaim{9.10. Proposition} Exact sequence (67) splits.\endproclaim 
\demo{Proof} Denote by $C$ the chain complex
$C_\ast(\tilde X)$ and by $B_p$, $Z_p$ and $H_p$ the boundaries, cycles,
and the homology of this chain complex, respectively. Since $\C[\pi]$
has homological dimension one, any submodule of a free module is free.
Thus, we have the free resolution
$$0\to B_p\to Z_p\to H_p\to 0\tag68$$
and also, there is a splitting
$$C_p\ =\ Z_p\ \oplus\ X_p$$
where $X_p$ is a free submodule such that the boundary homomorphism maps it
isomorphically onto $B_{p-1}$. Thus, we obtain
$$\Ker[M\otimes_\pi C_p\ \to\ M\otimes_\pi C_{p-1}]\ =
M\otimes_\pi Z_p\quad \oplus \quad 
\Ker[M\otimes_\pi X_p\to M\otimes_\pi Z_{p-1}].\tag69$$
Clearly, $\Ker[M\otimes_\pi X_p\to M\otimes_\pi Z_{p-1}]$ is isomorphic to
$\TOR_1^\pi(M,H_{p-1}(\tilde X))$ according to the definition 9.1.
By example 2.9, the extended $L^2$ homology $\H_p(X,M)$ is represented by the 
following morphisms of $\E(\A)$:
$$
\aligned
&(M\otimes_\pi C_{p+1}\to \Ker[M\otimes_\pi C_p\to M\otimes_\pi C_{p-1}])\ 
\simeq\ \\
&(M\otimes_\pi C_{p+1}\to M\otimes_\pi Z_p)\ \oplus\ 
\TOR_1^\pi(M,H_{p-1}(\tilde X))\ \simeq \\
&\TOR_0^\pi(M,H_p(\tilde X)) \oplus \TOR_1^\pi(M,H_{p-1}(\tilde X)).
\endaligned
$$                                      
This completes the proof. $\square$
\enddemo
\proclaim{9.11. Corollary} If the fundamental group $\pi=\pi_1(X)$ is free,
then for any $p$ the torsion part of the extended $L^2$ homology in dimension
$p$ coincides with the torsion part of $\TOR_0^\pi(M,H_p(\tilde X))$ and, 
in particular,
it depends only on the $\C[\pi]$ homology module $H_p(\tilde X)$ 
of the universal covering $\tilde X$. Thus, the Novikov-Shubin  
invariant of $\H_p(X,M)$ depends only on $H_p(\tilde X)$. $\square$
\endproclaim

Thus, in the case of the free fundamental group the Novikov-Shubin invariants
depend only on the homology modules of the universal covering $H_p(\tilde X)$.  

An example, computed by the author jointly with J.Hillman, shows
that for $\pi=\Z^2$ the Novikov-Shubin invariants are not functions of homology
of the universal covering (unlike the case of the free group) and depend also 
on the $k$-invariants.

A typical application of the spectral sequence of Theorem 9.7 consists in
getting estimates from above for the capacity of the extended $L^2$ homology.
The following statement yields an example. 

\proclaim{9.12. Theorem} Suppose that under the conditions of Theorem 9.7
it is known that all Hilbertian modules 
$$E^2_{p,q}\ =\ \TOR^\pi_q(M,H_p(\tilde X))\tag70$$
are torsion and have capacity less or equal than some $\nu\in [0,\infty)$.
Then the extended $L^2$ homology $\H_\ast(X,M)$ is also torsion and its
capacity is less or equal than $\nu$.
\endproclaim
\demo{Proof} It follows from Theorem 9.7 and Corollary 4.10.
$\square$\enddemo

\end